\documentclass[11pt]{amsart}
\usepackage{amssymb}
\usepackage{amsmath}	
\usepackage{ bm}

\usepackage[utopia]{mathdesign}

 \calclayout

\def\ignore#1{\relax}

\usepackage{graphicx}
\usepackage{epstopdf}
\newcommand\inlinegraphic[2][{scale=1.0}]{\begin{array}{c} \includegraphics[#1]{./EPS/#2}\end{array}}

\usepackage[plainpages=false, pdfpagelabels,  colorlinks=true, pdfstartview=FitV, linkcolor=blue, citecolor=blue, urlcolor=blue]{hyperref}

\numberwithin{equation}{section}
\numberwithin{figure}{section}

\def\la{\lambda}
\def\La{\Lambda}

\def\B{\mathfrak{B}}
\def\S{\mathfrak S}

\def\Ind{{\rm Ind}}
\def\Res{{\rm Res}}
\def\inv{^{-1}}

\def\uA{\breve A^\la}

\def\1{\bm 1}

\def\C{\mathcal{C}}

\def\spn{{\rm span}}

\def\br{B}
\def\mathbold{\bm}
\def\p #1{ \bm {#1}}
\def\pbar #1{\overline{\p {#1}}}

\def\inv{^{-1}}

\def\Z{{\mathbb Z}}

\def\Q{{\mathbb Q}}
\def\C{{\mathbb C}}
\def\la{\lambda}

\def\spp #1{^{(#1)}}

\def\ubold{{\bm u}}

\def\lambdabold{{\bm \la}} 
\def\labold{\lambdabold}
\def\mubold{{\bm \mu}}
\def\nubold{{\bm \nu}}

\def\qbold{{\bm q}}

\def\rhobold{{\bm \rho}}

\def\deltabold{{\mathbold \delta}}

\def\ahec #1{\widehat{H}_{#1}}  
\def\hec #1{H_{#1}}
\def\bmw #1{W_{#1}}

\def\tl{A}

\def\cont #1 #2 #3{C_{#1}^{#2}(#3)}

\theoremstyle{plain}
\newtheorem{theorem}{Theorem}[section]

\theoremstyle{plain}
\newtheorem{proposition}[theorem]{Proposition}

\theoremstyle{plain}
\newtheorem{corollary}[theorem]{Corollary}

\theoremstyle{plain}
\newtheorem{lemma}[theorem]{Lemma}

\theoremstyle{definition}
\newtheorem{definition}[theorem]{Definition}
\theoremstyle{definition}
\newtheorem{example}[theorem]{Example}

\theoremstyle{definition}
\newtheorem{remark}[theorem]{Remark}

\theoremstyle{definition}

\theoremstyle{definition}
\newtheorem{notation}[theorem]{Notation}

\theoremstyle{claim}

\theoremstyle{plain}
\newtheorem{assumption}[theorem]{Assumption}

    \def\mf{\mathfrak}

\title{On cellular algebras with Jucys Murphy elements}

\author{Frederick M. Goodman}
\address{Department of Mathematics\\ University of Iowa\\ Iowa
City, Iowa}
\email{ goodman@math.uiowa.edu} 

\author{John Graber}
\address{Department of Mathematics\\ University of Iowa\\ Iowa
City, Iowa}
\email{ jgraber@math.uiowa.edu}

\thanks{We thank the referee for very helpful criticism, which led to substantial improvements in the exposition.}


\begin{document}
 \begin{abstract}    We study analogues of Jucys-Murphy elements in cellular algebras arising from repeated Jones basic constructions.  Examples include Brauer and BMW algebras and their cyclotomic analogues. 
\end{abstract}
 \maketitle
\medskip

\setcounter{tocdepth}{1}
\tableofcontents

\section{Introduction}

We recently developed a framework for proving cellularity of a tower of algebras $(A_n)_{n\ge 0}$   that is  obtained from another tower of cellular  algebras $(Q_n)_{n\ge 0}$   by repeated Jones basic constructions
~\cite{GG1}.  
A key idea in this work is that of a tower of  algebras with coherent cellular structures;  coherence means that the cellular structures are well--behaved with respect to induction and restriction.
This paper continues our work on the themes of ~\cite{GG1};   here we
refine our framework by taking into account the role played by Jucys--Murphy elements. 

Before restricting to the  setting of ~\cite{GG1}, we first obtain some simple general results regarding coherent towers.  We show the existence of special cellular bases, called {\em path bases} which are distinguished by a restriction rule for the action of subalgebras on the basis elements.  We then give an axiomatization of Jucys--Murphy elements in coherent towers;  our assumptions imply that the JM elements act via triangular matrices on a path basis, as in Andrew Mathas's  axiomatization  \cite{mathas-seminormal} of cellular algebras with Jucys--Murphy elements.

\ignore{
 At the same time, we give a new version of  Andrew Mathas's  axiomatization  \cite{mathas-seminormal} of cellular algebras with Jucys--Murphy elements,  taking into account coherence of a sequence of such algebras.  While Mathas posits the triangularity property of the action of the Jucys--Murphy elements on the cellular basis, we derive this property from simpler assumptions.  
 }

Passing to the setting of ~\cite{GG1},  we use the general results mentioned above to 
give conditions which allow lifting Jucys--Murphy elements from $Q_n$ to $A_n$.  
Examples of algebras covered by this theory are Jones--Temperley--Lieb algebras, Brauer algebras, BMW algebras,
 and their cyclotomic analogues. Our method yields an easy and uniform proof of the triangularity property of the action of the Jucys--Murphy elements  in these examples, recovering theorems of Enyang ~\cite{Enyang2} and of Rui and Si ~\cite{rui-2008b} and ~\cite{rui-si-degnerate}.  
 
 \subsection{Antecedents and motivations}   Aside from our own paper ~\cite{GG1}, the most immediate antecedent and inspiration for this work was the paper of Andrew Mathas  ~\cite{mathas-seminormal} on Jucys--Murphy elements in cellular algebras.  As ~\cite{GG1} is about lifting cellular structures from Hecke--like algebras to BMW--like algebras, our intention was to find a way to lift Jucys--Murphy elements as well.   In order to do this, we needed a new axiomatization of Jucys--Murphy elements well adapted to the context of coherent towers.  The axiomatization that we propose does not replace that of Mathas, but compliments it;  a set of Jucys--Murphy elements in our sense is also a set of Jucys--Murphy elements in the sense of Mathas.

The Jucys--Murphy elements in $\C S_n$  were introduced by Murphy   ~\cite{murphy-seminormal-1981}in order to give a new  construction of Young's seminormal representations.  
The Jucys--Murphy elements of  $\C S_n$ generate the ``Gelfand--Zeitlin algebra" for the sequence  $(\C \mathfrak S_k)_{k \le n}$,  see section \ref{subsection: GZ algebras}; this is 
 a maximal abelian subalgebra of $\C \mathfrak S_n$ containing a canonical family of mutually orthogonal minimal idempotents $F_\mathfrak t$  indexed by Young tableaux of size $n$.     The seminormal basis of a simple module $\Delta^\lambda$  is obtained by a particular choice of one non--zero vector in the range of each $F_\mathfrak t$ for $\mathfrak t$ of shape $\la$.  This interpretation of the seminormal representations has been stressed by Ram ~\cite{ram-seminormal} and by  Okounkov  and Vershik  ~\cite{okounkov-vershik-selecta, okounkov-vershik-JMS-NY}. 

The JM elements in our theory duplicate this behavior; in a ``generic" setting, when the JM elements satisfy the separating condition of Mathas (and the algebras are in particular semisimple), 
 our JM elements generate the Gelfand--Zeitlin subalgebra for $(A_k)_{k \le n}$  for each $n$, see Proposition \ref{proposition: JM elements generate GZ algebra}.

\section{Preliminaries}  \label{section: preliminaries }

\subsection{Algebras with involution, and their bimodules}  \label{subsection: algebras with involution}
Let $R$ be a commutative ring with identity.
Recall that an involution $i$ on an $R$--algebra $A$ is   an $R$--linear algebra 
anti--automorphism of $A$ with \def\id{{\rm id}} $i^2 = \id_{A}$.  If $A$  and $B$ are 
$R$--algebras and $\Delta$ is an $A$--$B$ bimodule, then we define  a $B$--$A$ bimodule
$i(\Delta)$ as follows. As an $R$--module, $i(\Delta)$ is a copy of $\Delta$ with elements marked with the symbol $i$. The $B$--$A$ bimodule structure is defined by
$bi(x)a = i(i(a)xi(b))$.  Then $i$ is a functor  from the category of $A$--$B$ bimodules  to the category of $B$--$A$ bimodules.   By the same token,  we have a functor $i$ from the category of $B$--$A$ bimodules to the category of $A$--$B$ bimodules, and for an  $A$--$B$ bimodule
 $\Delta$,  we can identify $i\circ i(\Delta)$ with  $\Delta$. 
 
 Suppose that $A$, $B$, and $C$  are $R$--algebras with involutions  $i_A$,  $i_B$,  and $i_C$.  Let
$_B P_A$ and $_A Q_C$ be bimodules.   Then
$$
i(P \otimes_A Q) \cong i(Q) \otimes_A i(P),
$$
as $C$--$B$--bimodules.  Note that if we identify $i(P \otimes_A Q)$ with $i(Q) \otimes_A i(P)$, then we have the 
formula $i(p \otimes q) = i(q) \otimes i(p)$.  In particular,  let $M$ be a $B$--$A$--bimodule,  and identify
$i\circ i(M)$ with $M$, and $i(M\otimes_A i(M))$  with $i\circ i(M) \otimes_A i(M) = M\otimes_A i(M)$.
Then we have the formula $i(x \otimes i(y)) = y \otimes i(x)$.  

\subsection{Cellularity}
The definition of cellularity that we use is slightly weaker than the original definition of Graham and Lehrer in ~\cite{Graham-Lehrer-cellular}, see Remark \ref{remark:  on definition of cellularity}.

\begin{definition}  \label{gl cell}  Let $R$ be an integral domain and $A$ a unital $R$--algebra.  A {\em cell datum} for $A$ consists of  an algebra involution $i$ of $A$; a finite partially ordered set $(\Lambda, \ge)$ and 
for each $\la \in \Lambda$  a finite set  $\mathcal T(\lambda)$;  and   a subset $
\mathcal C = \{ c_{\mf s, \mf t}^\la :  \la \in \Lambda \text{ and } \mf s,\mf t \in \mathcal T(\la)\} \subseteq A$; 
with the following properties:
\begin{enumerate}
\item  $\mathcal C$ is an $R$--basis of $A$.
\item   \label{mult rule} For each $\la \in \Lambda$,  let $\breve A^\la$  be the span of the  $c_{\mf s,\mf  t}^\mu$  with
$\mu > \la$.   Given $\la \in \Lambda$,  $\mf s \in \mathcal T(\la)$, and $a \in A$,   there exist coefficients 
$r_{\mf v}^{\mf s}( a) \in R$ such that for all $\mf t \in \mathcal T(\la)$:
$$
a c_{\mf s, \mf t}^\la  \equiv \sum_{\mf v} r_{\mf v}^{\mf s}(a)  c_{\mf v,\mf  t}^\la  \mod  \breve A^\la.
$$
\item  $i(c_{\mf s,\mf  t}^\la) \equiv c_{\mf t, \mf s}^\la   \mod  \breve A^\la$ for all $\la\in \Lambda$ and $\mf s,\mf t \in \mathcal T(\lambda)$.

\end{enumerate}
$A$ is said to be a {\em cellular algebra} if it has a  cell datum.  
\end{definition}

For brevity,  we will write that  $(\mathcal C, \La)$ is a cellular basis of $A$.

\vbox{
\begin{remark} \mbox{} \label{remark:  on definition of cellularity}
\begin{enumerate}
\item  The original definition in  ~\cite{Graham-Lehrer-cellular} requires that $i(c_{\mf s,\mf t}^\la) = c_{\mf t, \mf s}^\la $ for all $\la,\mf s, \mf t$.  However, one can check that all of  \cite{Graham-Lehrer-cellular} remains valid with our weaker axiom.
\item  In case $2 \in R$ is invertible, one can check that our definition is equivalent to the original; ; see ~\cite{GG1}, Remark 2.4.
\item  One reason for using the weaker definition is that it allows a more graceful treatment of extensions of cellular algebras; see ~\cite{GG1}, Remark 2.6.  Another reason is that it becomes trivial to lift bases of cell modules to cellular bases of the algebra; see Lemma \ref{lemma: on globalizing bases of cell modules} and Remark \ref{remark: on globalizing bases} below. 
\end{enumerate}
\end{remark}
}

We recall some basic structures related to cellularity, see ~\cite{Graham-Lehrer-cellular}.
Given $\la\in\La$,   let $A^\la$ denote the span of the $c_{\mf s,\mf t}^{\mu}$ with $\mu \geq \la$.  It follows that both $A^\la$ and $\breve A^\la$ (defined above) are $i$--invariant two sided ideals of $A$.
\ignore{If $\mf t \in \mathcal T(\la)$, define $C_{\mf t}^\la$ to be the $R$-submodule of $A^\la/\uA$ with basis $\{ c_{\mf s,\mf t}^\la + \uA : \mf s \in \mathcal T(\la) \}$.  Then $C_{\mf t}^\la$ is a left $A$-module by Definition \ref{gl cell} (\ref{mult rule}).  Furthermore, the action of $A$ on $C_{\mf t}^\la$ is  independent of $\mf t$, i.e $C_{\mf u}^{\la}\cong C_{\mf t}^{\la}$  for any $\mf u,\mf t \in \mathcal T(\la)$.  }  The {\em left  cell module} $\Delta^\la$ 
is defined as follows: as   an $R$--module, $\Delta^\la$  is free with basis indexed by $\mathcal T(\la)$, say  $\{c_{\mf s}^\la$ : $\mf s \in \mathcal T(\la)\}$;   for each $a \in A$, the action of $a$ on $\Delta^\la$ is defined by  $ ac_{\mf s}^\la=\sum_{\mf v} r_{\mf v}^{\mf s}(a)  c_{\mf v}^\la$ where $r_{\mf v}^{\mf s}(a)$ is as  in Definition \ref{gl cell} (\ref{mult rule}).  
 
 For each $\la \in \Lambda$, we have an $A$--$A$--bimodule isomorphism $\alpha^\la : A^\la/\uA \rightarrow \Delta^\la \otimes_R i(\Delta^\la)$ determined by $\alpha^\la(c_{\mf s,\mf t}^{\la}+\uA)=c_{\mf s}^\la \otimes  i(c_{\mf t}^\la)$ satisfying   
 $i \circ \alpha^\la = \alpha^\la \circ i$,  using the remarks at the end of Section  \ref{subsection: algebras with involution} and point (2) of Definition \ref{gl cell}.   
 
 \subsection{Globalizing bases of cell modules}
 
 A given cellular algebra can have many cellular bases yielding the same cell modules and ideals $A^\la$.  The following lemma shows that an arbitrary collection of bases of the cell modules can be globalized to a cellular basis of the algebra.

 \begin{lemma}  \label{lemma: on globalizing bases of cell modules}
 
 Let $A$ be a cellular algebra, with cell datum denoted as above.  For each $\la \in \Lambda$, fix an $A$--$A$--bimodule isomorphism $\alpha^\la : A^\la/\uA \rightarrow \Delta^\la \otimes_R i(\Delta^\la)$ satisfying  $i \circ \alpha^\la = \alpha^\la \circ i$.     
 For each  $\la \in \Lambda$, let $\mathcal B^\la = \{b_\mf s^\la : \mf s \in \mathcal T(\la)\}$  be an arbitrary $R$--basis of 
 $\Delta^\la$.    For $\mf s, \mf t \in \mathcal T(\la)$,  let $b_{\mf s, \mf t}^\la$  be an arbitrary lifting of
 $(\alpha^\la)\inv(b_\mf s^\la \otimes b_\mf t^\la)$  to $A^\la$.  Then 
 $$
 \mathcal B = \{b_{\mf s, \mf t}^\la : \la \in \Lambda; \  \mf s, \mf t \in \mathcal T(\la)\}
 $$
 is a cellular basis of $A$.
 \end{lemma}
 
 \begin{proof}  It is easy to check that for each $\la \in \Lambda$,   
$ \{b_{\mf s, \mf t}^\mu : \mu \ge \la;  \ \mf s, \mf t \in \mathcal T(\mu)\}$  spans $A^\la$. 
In fact, if $\la$ is maximal in $\Lambda$, then   $A^\la \cong A^\la/\breve A^\la$,  and
$\{b_{\mf s, \mf t}^\la : \ \mf s, \mf t \in \mathcal T(\la)\}$ is a basis of $A^\la$.     Now fix $\la \in \Lambda$ and assume inductively that  for each $\la' > \la$,  
$ \{b_{\mf s, \mf t}^\mu : \mu \ge \la';  \ \mf s, \mf t \in \mathcal T(\mu)\}$  spans $A^{\la'}$.   This means that  
$ \{b_{\mf s, \mf t}^\mu : \mu > \la;  \ \mf s, \mf t \in \mathcal T(\mu)\}$  spans $\breve A^{\la}$. 
Now if $x \in A^\la$,  then $x \in  \spn  \{b_{\mf s, \mf t}^\la :  \mf s, \mf t \in \mathcal T(\la)\} + \breve A^\la$ and hence  $x \in \spn  \{b_{\mf s, \mf t}^\mu : \mu \ge \la;  \ \mf s, \mf t \in \mathcal T(\mu)\}$. 

Now it follows that $ \{b_{\mf s, \mf t}^\la : \la \in \Lambda; \  \mf s, \mf t \in \mathcal T(\la)\}$  spans
$A$.   Since $R$ is an integral domain, and this set has the same cardinality as the basis $\mathcal C$ of $A$,   it follows that the set is an $R$--basis of $A$.    Moreover, we have checked that $\breve A^\la$
(defined in terms of the original basis $\mathcal C$)  is the span of  the $b_{\mf s \mf t}^\mu$  with
$\mu > \la$.    Properties (2) and (3) of Definition \ref{gl cell}  (with $\mathcal C$  replaced by $\mathcal B$)   follow from the properties of the maps $\alpha^\la$.  
 \end{proof}
 
 \begin{remark}  \label{remark: on globalizing bases}
 Note that the proof only yields the weaker property (3) of Definition \ref{gl cell}   rather than the stronger requirement
$i(b_{\mf s,\mf t}^\la) = b_{\mf t, \mf s}^\la $ of ~\cite{Graham-Lehrer-cellular},  so this lemma would not be valid with the original definition of ~\cite{Graham-Lehrer-cellular}.  
\end{remark}

\begin{definition} \label{definition: globalization of bases}
If $\mathcal B^\la = \{b_\mf s^\la : \mf s \in \mathcal T(\la)\}$,  $\la \in \Lambda$ is a family of bases of the cell modules $\Delta^\la$, and
$\mathcal B =  \{b_{\mf s, \mf t}^\la : \la \in \Lambda; \  \mf s, \mf t \in \mathcal T(\la)\}$, is a cellular basis of $A$ such that $\alpha^\la(b_{\mf s, \mf t}^\la + \breve A^\la) = b_\mf s^\la \otimes b_\mf t^\la$   for each $\la, \mf s, \mf t$,  then we call $\mathcal B$ a {\em globalization} of the family of bases $\mathcal B^\la$,  $\la \in \Lambda$.  
\end{definition}

\subsection{Coherent towers of cellular algebras}
In ~\cite{GG1}, we defined a coherent tower of cellular algebras as follows:

\begin{definition} \label{definition: coherent tower}
Let $A_0 \subseteq A_1 \subseteq A_2 \subseteq \cdots$ be an increasing sequence of cellular algebras over an integral domain $R$.   Let $\Lambda_n$ denote the partially ordered set in the cell datum for $A_n$.   We say that $(A_n)_{n \ge 0}$ is a {\em  coherent tower of cellular algebras} if the following conditions are satisfied:
\begin{enumerate}
\item  The involutions are consistent; that is,  the involution on $A_{n+1}$,  restricted to $A_n$, agrees with the involution on $A_n$.
\item  For each $n\ge 0$ and for each $\la \in \Lambda_n$, the induced module $\Ind_{A_n}^{A_{n+1}} (\Delta^\la)$
has a filtration by cell modules of $A_{n+1}$. That is, there is a filtration
$$
\Ind_{A_n}^{A_{n+1}} (\Delta^\la) = M_t \supseteq M_{t-1} \supseteq \cdots \supseteq M_0 = (0)
$$
such that for each $j\ge1$,  there is a $\mu_j \in \Lambda_{n+1}$  with $M_j/M_{j-1} \cong \Delta^{\mu_j}$.
\item  For each $n\ge 0$ and for each $\mu \in \Lambda_{n+1}$, the restriction  $\Res_{A_n}^{A_{n+1}} (\Delta^\mu)$
has a filtration by cell modules of $A_{n}$. That is, there is a filtration
$$
\Res_{A_n}^{A_{n+1}} (\Delta^\mu) = N_s \supseteq N_{s-1} \supseteq \cdots \supseteq N_0 = (0)
$$
such that for each $i\ge1$, there is a $\la_i \in \Lambda_{n}$  with $N_j/N_{j-1} \cong \Delta^{\la_i}$.

\end{enumerate}
\end{definition}

The modification of the definition for a {\em finite} tower of cellular algebras is obvious. We call a filtration as in (2) and (3) a {\em cell filtration}.
In the examples of interest to us, we will also have {\em uniqueness of the multiplicities} of the cell modules appearing as subquotients of the cell filtrations, and {\em Frobenius reciprocity} connecting the multiplicities in the two types of filtrations.  We did not include uniqueness of multiplicities and Frobenius reciprocity as requirements in the definition, as they will follow from additional assumptions that we will impose later.

We introduce a stronger notion of coherence:
\begin{definition}  \label{definition: strongly coherent tower}
Say that a coherent tower of cellular algebras $(A_n)_{n \ge 0}$ is  {\em  strongly coherent}  if 
$A_0 \cong R$ and 
in the cell filtrations (2) and (3) in Definition \ref{definition: coherent tower}, we have
$$
\mu_{t} < \mu_{t-1} < \cdots < \mu_{1}$$  in the partially ordered set $\La_{n+1}$,   and
$$\la_{s} < \la_{s-1} < \cdots < \la_{1}$$  in the partially ordered set $\La_{n-1}$.  
\end{definition}

\subsection{Inclusions of split semisimple algebras and branching diagrams}

 A finite dimensional split semisimple algebra over a field $F$ is one which is isomorphic to a finite direct sum of full matrix algebras over $F$.
 
Suppose
$A \subseteq B$ are finite dimensional split semisimple algebras over  $F$ (with the same identity element).  Let   $A(i)$,  $i \in I$, be the minimal ideals of $A$  and  $B(j)$,  $j \in J$,   the minimal ideals of $B$.  
We associate a  $J \times I$   {\em inclusion matrix}
$\Omega$ to the inclusion $A \subseteq B$, as follows.  Let $W_j$ be a simple $B(j)$--module.
Then $W_j$ becomes an $A$--module via the inclusion,  and $\Omega(j, i)$ is defined to be the multiplicity of  a simple $A(i)$--module 
 in the decomposition of $W_j$ as an $A$--module.

It is convenient to encode an inclusion matrix   by a bipartite graph, called the {\em branching diagram};  the branching diagram has vertices labeled by $I$  arranged on one horizontal line,  vertices labeled by  $J$  arranged along a second (higher) horizontal line,    and $\Omega(j, i)$ edges connecting
$j \in J$ to $i \in I$.

If $A_0 \subseteq A_1 \subseteq A_2 \subseteq   \cdots$ is a (finite or infinite) sequence of inclusions of finite dimensional split semisimple algebras over $F$,  then the branching diagram for the sequence is obtained by stacking the branching diagrams for each inclusion,  with the 
upper vertices of the diagram for $A_i \subseteq A_{i+1}$  being identified with the lower vertices of the diagram for $A_{i+1} \subseteq A_{i+2}$.
For two vertices $\la$ on level $\ell$ of a branching diagram and  $\mu$ on level $\ell + 1$, write 
$\la \nearrow \mu$  if $\la$ and $\mu$ are connected by an edge.

\begin{notation}
Let $R$ be an integral domain with field of fractions $F$.  Let $A$ be a cellular algebra over $R$ and $\Delta$ an $A$--module.  Write $A^F$  for $A \otimes_R F$  and
$\Delta^F$ for $\Delta \otimes_R F$.  
\end{notation}
 
 \begin{lemma}[\cite{GG1}, Lemma 2.20] \label{lemma: multiplicities in cell filtrations}
Let $R$ be an integral domain with field of fractions $F$.
 Suppose that $(A_n)_{n \ge 0}$ is a coherent tower of cellular algebras over  $R$ and that
 $A_n^F$ is split semisimple for all $n$.  Let $\Lambda_n$ denote the partially ordered set in the cell datum for
 $A_n$. 
  Then
 \begin{enumerate}
 \item 
 $\{(\Delta^\la)^F : \la \in \Lambda_n\}$ is a complete family of simple $A_n^F$--modules.
 \item Let $[\omega(\mu, \la)]_{\mu \in \Lambda_{n+1}, \,  \la \in \Lambda_n}$ denote the inclusion matrix for
 $A_n^F \subseteq A_{n+1}^F$.   Then for any $\la \in \Lambda_n$ and $\mu \in \Lambda_{n+1}$, 
  and any cell filtration of $\Res_{A_n}^{A_{n+1}}(\Delta^\mu)$,  the number of subquotients of the filtration isomorphic to  $\Delta^\la$ is $\omega(\mu, \la)$.
  
  \item  Likewise, for any $\la \in \Lambda_n$ and $\mu \in \Lambda_{n+1}$, 
 and any cell filtration of $\Ind_{A_n}^{A_{n+1}}(\Delta^\la)$,  the number of subquotients of the filtration isomorphic to  $\Delta^\mu$ is $\omega(\mu, \la)$.
 \end{enumerate}
\end{lemma}

\begin{corollary}  Under the hypotheses of Lemma \ref{lemma: multiplicities in cell filtrations}, the multiplicity of a cell module as a subquotient of a cell filtration of $\Res_{A_n}^{A_{n+1}}(\Delta^\mu)$
or of $\Ind_{A_n}^{A_{n+1}}(\Delta^\la)$ is independent of the choice of the cell filtration.  Moreover,
The multiplicity of $\Delta^\la$ in $\Res_{A_n}^{A_{n+1}}(\Delta^\mu)$ equals the multiplicity of 
$\Delta^\mu$ in $\Ind_{A_n}^{A_{n+1}}(\Delta^\la)$.
\end{corollary}

\begin{definition} A tower of split semisimple algebras $(A_n)_{n \ge 0}$   over a field $F$ is {\em multiplicity free} if all entries in the inclusion matrices are $0$ or $1$ and $A_0 \cong F$.    Equivalently,  there are no multiple edges in the branching diagram of the tower, and there is a unique vertex (denoted $\emptyset$) at level $0$.   We will also say that the branching diagram is multiplicity free. 
\end{definition}

\begin{corollary}   Under the hypotheses of Lemma \ref{lemma: multiplicities in cell filtrations},  if
$(A_n)_{n\ge0}$ is strongly coherent, then $(A_n^F)_{n\ge0}$ is a multiplicity free tower of split semisimple algebras.
\end{corollary}

\begin{example}  \label{example: hecke algebras 1}
 Fix an integral domain $S$ and an invertible $q \in S$.   The Hecke algebra 
$H_{n}(q) = H_{n, S}(q)$  is the associative, unital $S$--algebra with generators $T_{j}$ for
$1 \le j \le n-1$,  satisfying the braid relations and the quadratic relation 
$(T_{j} -q)(T_{j} +1) = 0$ for all $j$.    $H_n(q)$ has an algebra involution $x \mapsto x^*$ uniquely determined by $(T_j)^* = T_j$.    $H_n(q)$  has a cellular basis due to Murphy
 ~\cite{murphy-hecke95}
$$
\{m_{\mathfrak s, \mathfrak t}^\la :  \la \in Y_n; \    \mathfrak s, \mathfrak t \in \mathcal T(\la)\},
$$
where $Y_n$ is the partially ordered set of all Young diagrams of size $n$, with dominance order $\unrhd$, and $\mathcal T(\la)$ is the set of all standard Young tableaux of shape $\la$.
By results of Murphy ~\cite{murphy-hecke95}, Dipper and James ~\cite{dipper-james1, dipper-james2}, and Jost ~\cite{jost},   the sequence of Hecke algebras
$(H_{n, S}(q))_{n \ge 0}$  is  strongly coherent.

The generic ground ring for the Hecke algebras is $R = \Z[\qbold, \qbold\inv]$,  where $\qbold$ is an indeterminant over $\Z$;  the Hecke algebra $H_{n, S}(q)$ over any $S$ is a specialization of $H_{n, R}(\qbold)$.  If $F = \Q(\qbold)$  denotes the field of fractions of $R$, then $H_{n, F}(\qbold)$  is split semisimple for all $n$  and the branching diagram for the tower of Hecke algebras $(H_{n, F}(\qbold))_{n \ge 0}$  is Young's lattice $\mathcal Y$, which is multiplicity free.  
\end{example}

\subsection{Remark on the role of generic ground rings}  In the examples of interest to us (Hecke algebras, BMW algebras, etc.)  there is a generic  ground ring $R$ with the properties that:
\begin{enumerate}
\item  $R$ is an integral domain and the algebras $A_n^F$ over the field of fractions of $R$ are split semisimple, and
\item  the algebras over any ground ring $S$ are specializations of those over $R$,  $A_n^S = A_n^R \otimes_R S$.   
\end{enumerate}

Certain properties of the algebras over the generic ground ring $R$ carry over to any specialization.  For example, if the algebras over $R$ are cellular, so are all of the specializations.   For another example, in the next section, we show the existence of certain bases, called path bases, in strongly coherent towers of cellular algebras over an integral domain $R$, assuming the algebras over the field of fractions of $R$ are semisimple.   This hypothesis would apply to the generic ground ring in our examples.  But then the path bases in the cell modules over $R$ can be specialized to cell modules over any ground ring $S$.

\subsection{Path bases in strongly coherent towers} 
\label{subsection: bases in strongly coherent towers}
 In this section, we discuss {\em path bases} in  strongly coherent towers of cellular algebras.

\begin{assumption}
In  Section \ref{subsection: bases in strongly coherent towers}, let  $R$ be an integral domain with field of fractions $F$, 
 $(A_n)_{n \ge 0}$ a strongly  coherent tower of cellular algebras over $R$, such that
 $A_n^F$ is semisimple for all $n$.  Let 
 $\B$  denote  the branching diagram of  $(A_n^F)_{n \ge 0}$ and   $\Lambda_n$  the partially ordered set in the cell datum for $A_n$. 
 \end{assumption}

\begin{definition}
A {\em path} on $\B$   from $\la \in \La_\ell$  to
$\mu \in \La_m$  ($\ell < m$)  is a sequence $(\la = \la\spp \ell, \la\spp {\ell + 1}, \dots,  \la\spp {m} = \mu)$  with $\la \spp i \nearrow \la \spp {i+1}$  for all $i$.  A path $\mathfrak s$  from $\la$ to 
$\mu$  and a path $\mathfrak t$  from $\mu$ to $\nu$  can be concatenated in the obvious way;  denote the concatenation $\mathfrak s \circ \mathfrak t$.   If
$\mathfrak t = (\emptyset = \la\spp 0, \la\spp 1, \dots, \la\spp n = \la)$ is a path from 
$\emptyset$ to $\la \in \La_n$, and $0 \le k < \ell \le n$,  write $\mf t(k) = \la\spp k$, 
$\mathfrak t_{[k, \ell]}$ for
the path  $(\la\spp k,  \dots, \la\spp \ell)$, and  write $\mathfrak t'$ for $\mathfrak t_{[0,n-1]}$.
\end{definition}

For $\la \in \La_n$, the rank of the cell module $\Delta^\la$ of $A_n$ is the same as the dimension of the simple $A_n^F$ module $(\Delta^\la)^F$, namely the number of paths
on $\B$ from $\emptyset$ to $\la$.  
It follows that we can assume without loss of generality that the index set  $\mathcal T(\la)$ in the cell datum for $A_n$ is 
equal to the set of paths
on $\B$ from $\emptyset$ to $\la$.   We set $\mathcal T(n) = \bigcup_{\la \in \La_n} \mathcal T(\la)$, the set of paths on $\B$ from $\emptyset$ to some $\la \in \La_n$.  

\begin{definition} \label{definition: orders on paths}
(Partial orders on the set of paths.)  We introduce two natural partial orders on $\mathcal T(n)$.  Let $\mathfrak s = (\la\spp{0} , \la\spp{1}, \dots, \la\spp{n} )$  and
$\mathfrak t = (\mu\spp{0} , \mu\spp{1}, \dots, \mu\spp{n})$  be two  paths with
$\la\spp i, \mu\spp i \in \Lambda_{i}$.    Say that $\mathfrak s$  precedes $\mathfrak t$ in {\em dominance order}  (denoted $\mathfrak s \unlhd \mathfrak t$)  if $\la\spp i \le \mu\spp i$ for all $i$ ($ 0 \le i \le n$).    Say that $\mathfrak s$  precedes $\mathfrak t$ in {\em reverse lexicographic  order}  (denoted $\mathfrak s \preceq \mathfrak t$)  if $\mathfrak s = \mathfrak t$, or if  for the last index $j$ such that $\la\spp j \ne \mu \spp j$, we have $\la\spp j < \mu \spp j$ in
$\Lambda_{j}$.  Similarly,  we can order the paths going from level $k$ to level $n$ on $\B$ 
by dominance or by reverse lexicographic order.  
\end{definition}  

\begin{example} Take $\B$ to be Young's lattice.  For a Young diagram $\la$, 
standard Young tableaux of shape $\la$  can be identified with paths on  $\B$  from the empty diagram to $\la$.   Dominance order on paths, as defined in Definition \ref{definition: orders on paths}, agrees with dominance order on standard tableaux as usually defined.    Reverse lexicographic order
  coincides with the ``last letter order," see for example   ~\cite{murphy-seminormal-1981}, page 288. 
  \end{example}
  
We will now construct  certain bases $\mathcal B^\la = \{b_\mathfrak s^\la : \mathfrak s \in \mathcal T(\la)\}$ of the cell modules $\Delta^\la$, $\la \in \cup_n \La_n$, each 
indexed by the set of paths $\mathcal T(\la)$,  by induction on $n$. For $\lambda \in \Lambda_0 $ or
$\lambda \in \Lambda_1 $, the cell module $\Delta^\lambda $ is free of rank one, and we choose any basis.  Suppose now that $n > 1$,  and a basis $\{b_\mathfrak s^\mu : \mathfrak s \in \mathcal T(\mu)\}$  for  $\Delta^\mu$  has been obtained for each $\mu \in \Lambda_k$ for $k \le n-1$.   Let $\la \in \Lambda_n$, and consider the filtration 
\begin{equation} \label{equation: filtration of restricted module 1}
\Res_{A_{n-1}}^{A_n}(\Delta^\la) =   N_s \supseteq N_{s-1} \supseteq \cdots \supseteq N_0 = (0),
\end{equation}
with $N_j/N_{j-1} \cong \Delta^{\mu_j}$  and $\mu_s < \mu_{s-1} < \cdots < \mu_1$.  
For each $j$, let $\{\bar b_\mathfrak s^{\mu_j} : \mathfrak s \in \mathcal T(\mu_j)\}$  be any lifting to $N_j$  of the basis
$\{b_\mathfrak s^{\mu_j}  : \mathfrak s \in \mathcal T(\mu_j)\}$ of $N_j/N_{j-1} \cong \Delta^{\mu_j}$.  Then $\cup_j \{\bar b_\mathfrak s^{\mu_j} : \mathfrak s \in \mathcal T(\mu_j)\}$ 
is a basis of $\Delta^\la$.    Note that $\mathfrak t \mapsto \mathfrak t'$ is a bijection from
$\mathcal T(\la)$ to $\cup_j  \mathcal T(\mu_j)$.   We define $b_\mathfrak t^\la$ to
be $\bar b_{\mathfrak t'}^{\mu_j}$  if  $t' \in \mathcal T(\mu_j)$, so our basis is now denoted by   
$\{b_\mathfrak s^\la : \mathfrak s \in \mathcal T(\la)\}$.
The bases $\mathcal B^\la = \{b_\mathfrak s^\la : \mathfrak s \in \mathcal T(\la)\}$ of the cell modules $\Delta^\la$  have the following property.

\begin{proposition} \label{proposition:  action of Hk on basis for strong coherence}

Fix $0\le k < n$,  $\la \in \La_n$, and $\mathfrak t \in \mathcal T(\la)$. Write
$\mu = \mathfrak t(k)$,  $\mathfrak t_1 = \mathfrak t_{[0,k]}$, and $\mathfrak t_2 = \mathfrak t_{[k,n]}$.   Let $x \in A_k$, and let $x b_{\mathfrak t_1}^\mu = \sum_{\mathfrak s} r(x; \mathfrak s, \mathfrak t_1)  b_{\mathfrak s}^\mu$.  Then
$$
x b_\mathfrak t^\la \equiv  \sum_{\mathfrak s} r(x; \mathfrak s, \mathfrak t_1)  b_{\mathfrak s \circ \mathfrak t_2}^\la,
$$ 
modulo  $\spn \{ b_\mathfrak v^\la   : \mathfrak v_{[k,n]} \succ \mathfrak t_{[k,n]}         \}$, where $\succ$ denotes reverse lexicographic order.
\end{proposition}

\begin{proof}  We prove this by induction on $n-k$.  Consider the case $n-k = 1$.  
Consider the filtration (\ref{equation: filtration of restricted module 1}). If $\mathfrak t' \in
\mathcal T(\mu_j)$,  then  by the construction of the basis $\{b_{\mathfrak t}^\la : \mathfrak t \in \mathcal T(\la)\}$, we have
$$
x b_\mathfrak t^\la \equiv  \sum_{\mathfrak s} r(x; \mathfrak s, \mathfrak t_1)  b_{\mathfrak s \circ \mathfrak t_2}^\la,
$$ 
modulo $N_{j-1}$,  while $N_{j-1}$ equals the $R$--span of $\{ b_\mathfrak v^\la   : \mathfrak v_{[n-1,n]} \succ \mathfrak t_{[n-1,n]}         \}$. 

 Now suppose that $n-k>1$, and $\mathfrak t' \in \mathcal T(\mu_j)$.  Then $x b_\mathfrak t^{\la} = x \bar b_{\mathfrak t'}^{\mu_j} $.  By  a suitable  induction hypothesis,
$$
x b_{\mathfrak t'}^{\mu_j} \equiv  \sum_{\mathfrak s} r(x; \mathfrak s, \mathfrak t_1)  b_{\mathfrak s \circ \mathfrak t_{[k, n-1]}}^{\mu_j},
$$ 
modulo the span of 
$\{ b_\mathfrak v^{\mu_j}   : \mathfrak v_{[k,n-1]} \succ \mathfrak t_{[k,n-1]}         \}$.  But then
$$
x b_\mathfrak t^\la \equiv  \sum_{\mathfrak s} r(x; \mathfrak s, \mathfrak t_1)  b_{\mathfrak s \circ \mathfrak t_2}^\la,
$$
modulo
$$
\spn \{ b_\mathfrak v^{\mu_j}   : \mathfrak v_{[k,n-1]} \succ \mathfrak t_{[k,n-1]}         \} + N_{j-1} = 
\spn \{ b_\mathfrak v^{\la}   : \mathfrak v_{[k,n]} \succ \mathfrak t_{[k,n]}\}  .
$$
\end{proof}

With the family of bases $\mathcal B^\la = \{b_\mathfrak s^\la : \mathfrak s \in \mathcal T(\la)\}$ of the cell modules $\Delta^\la$,  as above,    for each $n \ge 0$,  let $\mathcal B_n = \{b_{ \mf s, \mf t }^\la : \la \in \La_n, \ \mf s, \mf t \in \mathcal T(\la)\}$  be a cellular basis of 
$A_n$  globalizing the bases $\mathcal B^\la$,  $\la \in \La_n$;  see Lemma \ref{lemma: on globalizing bases of cell modules}
and Definition \ref{definition: globalization of bases}. 

The cellular bases $\mathcal B_n = \{b_{ \mf s, \mf t }^\la : \la \in \La_n, \ \mf s, \mf t \in \mathcal T(\la)\}$  have the property:

\begin{corollary}  Fix $0\le k < n$,  $\la \in \La_n$, and $\mf  t \in \mathcal T(\la)$. Write
$\mu = \mf  t(k)$,  $\mf  t_1 = \mf  t_{[0,k]}$, and $\mf  t_2 = \mf  t_{[k,n]}$.   Let $x \in A_k$, and let $x b_{\mf  t_1, \mf v}^\mu \equiv \sum_{\mf  s} r(x; \mf  s, \mf  t_1)  b_{\mf  s, \mf v}^\mu$ modulo $\breve A_k^\mu$ for all
$\mf v \in \mathcal T(\mu)$.   Then, for all $\mf v \in \mathcal T(\la)$, 
$$
x b_{\mf  t, \mf v} ^\la  \equiv  \sum_{\mf  s} r(x; \mf  s, \mf  t_1)  b_{\mf  s \circ \mf  t_2, \mf v}^\la,
$$ 
modulo  $\spn \{ b_{\mf  w, \mf v}^\la   : \mf  w_{[k,n]} \succ \mf  t_{[k,n]}         \} +  \breve A_n^\la$, 
\end{corollary}

\begin{definition} \label{definition: path bases}
 A family of bases $\mathcal B^\la$  of the cell modules $\Delta^\la$, $\la \in \bigcup_n \La_n$,  
having the property described in Proposition \ref{proposition:  action of Hk on basis for strong coherence}  will be called a {\em family of path bases} of the cell modules.   

A family of cellular bases $\mathcal B_n$   of $A_n$,  $n \ge 0$,  globalizing  a family of path bases $\mathcal B^\la$ of the cell modules     will also be called a  {\em family of path bases} of the cellular algebras.  

\end{definition}

\section{JM elements in coherent towers}  \label{subsection: preliminaries: JM}
\label{section: JM elements in coherent towers}
\begin{example}  \label{example: JM elements in Hecke algebra}
We recall the classical Jucys--Murphy elements  in the Hecke  algebra $H_{n}(q)$, and some of their properties.    
The (multiplicative) Jucys--Murphy elements in $H_n(q)$  are the elements
$\{L_1, \dots, L_n\}$  defined by $L_1 = 1$  and $L_{j+1} =q\inv T_j L_j T_j$   for $1 \le j \le n-1$.  The elements $L_k$  are mutually commuting;  in fact, $L_k \in H_k(q) \subseteq H_n(q)$ for $1 \le k \le n$,   and for $k \ge 2$,  $L_k$  commutes with $H_{k-1}$.   Symmetric polynomials in the $\{L_k\}$  are in the center of $H_n(q)$.    The  Jucys--Murphy elements act on the Murphy bases of the cell module $\Delta^\la$  as follows.  Let $\kappa(j, \mathfrak t) = c(j, \mathfrak t) - r(j, \mathfrak t)$, where $c(j, \mathfrak t)$ is the column of $j$ in the standard tableau $\mathfrak t$ and $ r(j, \mathfrak t)$ is the row of $j$ in $\mathfrak t$.  Then
 \begin{equation} \label{equation: triangular action of JM elements in Hecke example}
 L_j m_\mathfrak t^\la =  q^{\kappa(j, \mathfrak t)}  m_\mathfrak t^\la  + \sum_{\mathfrak s \rhd \mathfrak t} r_\mathfrak s  m_\mathfrak s^\la.
  \end{equation}
  For a cell $x$ in the Young diagram $\la$,  let $\kappa(x)$  denote its content, namely 
  the column of $x$ minus the row of $x$. 
  It follows from (\ref{equation: triangular action of JM elements in Hecke example})  that the product $p = \prod_{j = 1}^n L_j$  acts as a scalar $\alpha_\la = q^{\sum_{x \in \la}\kappa(x)}$  on the cell module $\Delta^\la$.  Namely,  if $\mathfrak t_0$ is the most dominant standard tableaux of shape $\la$ then $p m_{\mathfrak t_0}^\la = \alpha_\la m_{\mathfrak t_0}^\la $, by
  (\ref{equation: triangular action of JM elements in Hecke example}). But $p$ is central and
  $\Delta^\la$ is a cyclic module with generator $m_{\mathfrak t_0}^\la$.
  \end{example}

 Abstracting from the Hecke algebra example,  Mathas
  ~\cite{mathas-seminormal} defined a family of JM--elements in a cellular algebra as follows.
  
  \begin{definition}[\cite{mathas-seminormal}]  Let $A$ be a cellular algebra over $R$; let $\La$  denote the  partially ordered set in the cell datum for $A$, and, for each $\la \in \La$,  let $\{a_\mathfrak t^\la:  \mathfrak t \in \mathcal T(\la)\}$ denote the  basis  of the cell module $\Delta^\la$ (derived from the cellular basis of $A$.)    Suppose 
that for each $\la \in \La$,    the index set $\mathcal T(\la)$ is given a partial order $\succeq$. 

  A finite family of elements $\{L_j: 1 \le j \le M\}$ in $A$  is a {\em JM--family in the sense of Mathas} if the elements $L_j$  are mutually commuting and invariant under the involution of $A$, and,   for each $\la \in \La$,  
  there is a set of scalars $\{\kappa(j, \mathfrak t) : 1 \le j \le n,  \mathfrak t \in \mathcal T(\la)\}$  such that for $1 \le j \le n$ and  $\mathfrak t \in \mathcal T(\la)$,
  $$
  L_j a_\mathfrak t^\la = \kappa(j, \mathfrak t) a_\mathfrak t^\la + \sum_{\mathfrak s\, \succ \, \mathfrak t} r_\mathfrak s a_\mathfrak s^\la,
  $$
  for some $r_\mathfrak s \in R$, depending on $j$ and $\mathfrak t$.  In addition, the family $\{L_j\}$ is said to be {\em separating}  if 
 $\mathfrak t \mapsto (\kappa(j, \mathfrak t))_{1 \le j \le n}$ is injective on ${\mathcal T = \bigcup_{\la \in \La} \mathcal T(\la)}$.\footnote{Mathas' definition of separating is slightly weaker.}  
 \end{definition}
  
 We are going to introduce  a different abstraction of Jucys--Murphy elements that is appropriate for  strongly coherent towers of cellular algebras.  We will see that  our concept implies that of Mathas.    
 
  \begin{definition} \label{definition: multiplicative JM family}
   Let   $(A_n)_{n \ge 0}$ be a strongly coherent tower of cellular algebras over $R$. Let $\La_n$  denote the partially ordered set in the cell datum for $A_n$.
  
 A family of {\em invertible} elements $\{L_n:  n \ge 1\}$  is a {\em multiplicative JM--family}  if for all $n \ge 1$,
 \begin{enumerate}
  \item $L_n \in A_n$, $L_n$ is invariant under the involution of $A_n$,   and, for $n \ge 1$, 
  $L_n$ commutes with $A_{n-1}$.  In particular, the elements $L_j$ are mutually commuting.
  \item For each $n \ge 1$ and each $\la \in\La_n$,   there exists an invertible $\alpha(\la) \in R$  such that  the product $ L_1 \cdots L_n$   acts as the scalar   $\alpha(\la)$  on the cell module $\Delta^\la$.
    \end{enumerate}
  \end{definition}
  
  For convenience, we will set $\alpha(\emptyset) = 1$,  where $\emptyset$ is the unique element of $\Lambda_0$.

  \begin{definition}  \label{definition: additive JM family}
  An {\em additive JM--family}  is defined similarly,
  except that the elements $L_j$ are not required to be invertible and (2) is replaced by 
  
 \     (2$'$) \     For each $n \ge 1$ and each $\la \in\La_n$,   there exists  $d(\la) \in R$  such that  the sum $ L_1 +  \cdots + L_n$   act as the scalar   $d(\la)$  on the cell module $\Delta^\la$.
  \end{definition}
  
  For convenience, we will set $d(\emptyset) = 0$.
    
  \begin{assumption} \label{assumption: assumptions on coherent tower for section 3}
  For the remainder of  Section \ref{section: JM elements in coherent towers}, let  $R$ be an integral domain with field of fractions $F$, 
 $(A_n)_{n \ge 0}$ a strongly coherent tower  of cellular algebras over $R$, such that $A_n^F$ is semisimple for all $n$.  Let    $\B$  denote the branching diagram of  $(A_n^F)_{n \ge 0}$ and    $\Lambda_n$  the partially ordered set in the cell datum for $A_n$.  
Let $\mathcal B^\la = \{b_\mathfrak s^\la : \mathfrak s \in \mathcal T(\la)\}$ be a family of path bases of the cell modules $\Delta^\la$, $\la \in \bigcup_n \La_n$ (Definition \ref{definition: path bases}).     We employ the reverse lexicographic order $\preceq$ on paths (Definition
 \ref{definition: orders on paths}). 
 \end{assumption}

 \begin{proposition} \label{proposition: triangularity property of JM elements}
Suppose that $\{L_n: n \ge 0\}$  is a multiplicative JM--family for the strongly coherent tower $(A_n)_{n \ge 0}$.
\begin{enumerate}
\item
For $n \ge 1$  and $\la \in \La_n$, let 
 $\alpha(\la) \in R^\times$  be  such that $ L_1 \cdots L_n$  acts by the scalar $\alpha(\la)$ on the cell module $\Delta^\la$.  Then for all $n \ge 1$,  $\la \in \La_n$,   $\mathfrak t \in \mathcal T(\la)$, and $1 \le j \le n$,  we have
 \begin{equation} \label{equation: triangular action of Lj 1}
   L_j b_\mathfrak t^\la = \kappa(j, \mathfrak t)\  b_\mathfrak t^\la + \sum_{\mathfrak s\, \succ \, \mathfrak t} r_\mathfrak s b_\mathfrak s^\la,
 \end{equation}
  for some elements $r_\mathfrak s \in R$ (depending on $j$ and $\mathfrak t$), 
 with $\kappa(j, \mathfrak t) = \displaystyle\frac{\alpha(\mathfrak t(j))}{\alpha(\mathfrak t(j-1))}$.  
   \item For each $n \ge 1$,  $ L_1 \cdots L_n$ is in the center of $A_n$.  
   \end{enumerate}  
 \end{proposition}
 
 \begin{proof}  We prove (1)  by induction on $n$. For $n=1$,   the statement follows from 
 (2) of Definition \ref{definition: multiplicative JM family}.   Assume $n>1$ and adopt the appropriate induction hypothesis.  For $j<n$,  $\la \in \Lambda_n$, and $\mathfrak t \in \mathcal T(\lambda) $,  (\ref{equation: triangular action of Lj 1}) holds by the induction hypothesis and Proposition \ref{proposition:  action of Hk on basis for strong coherence}, while
 $$
 \begin{aligned}
L_n b_\mathfrak t^\la &= ( L_1 \cdots L_{n-1})\inv (L_1 \dots L_n)  b_\mathfrak t^\la\\
& =\alpha(\la) \ ( L_1 \cdots L_{n-1})\inv  b_\mathfrak t^\la\\
&= \alpha(\la) \alpha(\mathfrak t(n-1))\inv b_\mathfrak t^\la  + 
 \sum_{\mathfrak s\, \succ \, \mathfrak t} r_\mathfrak s b_\mathfrak s^\la,
 \\
\end{aligned}
 $$
  using point
 (2) of Definition \ref{definition: multiplicative JM family} and  
  Proposition \ref{proposition:  action of Hk on basis for strong coherence}.
  
 For all $x \in A_n$,  $x( L_1 \cdots L_n) = (L_1 \cdots L_n)x$ on each cell module. But the direct sum of all cell modules is faithful.   This proves (2).  
 \end{proof}
 
   The additive version  of the proposition is the following;  the proof is similar. 
Recall that Assumption \ref{assumption: assumptions on coherent tower for section 3} is still in force.    
    \begin{proposition} \label{proposition: triangularity property of  additive JM elements}
Suppose that $\{L_n: n \ge 0\}$  is an additive JM--family for the  tower $(A_n)_{n \ge 0}$.
\begin{enumerate}
\item
For $n \ge 1$  and $\la \in \La_n$, let 
 $d(\la) \in R$  be  such that $L_1 + \cdots + L_n$  acts by the scalar $d(\la)$ on the cell module $\Delta^\la$.  Then for all $n \ge 1$,  $\la \in \La_n$,   $\mathfrak t \in \mathcal T(\la)$, and $1 \le j \le n$,  we have
 \begin{equation} \label{equation: triangular action of Lj 2}
   L_j b_\mathfrak t^\la = \kappa(j, \mathfrak t) b_\mathfrak t^\la + \sum_{\mathfrak s\, \succ \, \mathfrak t} r_\mathfrak s b_\mathfrak s^\la,
 \end{equation}
  for some elements $r_\mathfrak s \in R$  (depending on $j$ and $\mathfrak t$),
 with $\kappa(j, \mathfrak t) = \displaystyle{\alpha(\mathfrak t(j))}- {\alpha(\mathfrak t(j-1))}$. 
   \item For each $n \ge 1$,  $ L_1 +\cdots + L_n$ is in the center of $A_n$.  
   \end{enumerate}  
 \end{proposition} 
 
 \begin{remark}  The techniques employed here give triangularity of the action of the JM elements only with respect to the reverse lexicographic order on paths, and not with respect to the dominance order.  Our techniques cannot recover the result on triangularity with respect to the dominance order for the Hecke algebras (see Example \ref{example: JM elements in Hecke algebra}).    
 \end{remark}
 
 \subsection{The separated case -- Gelfand--Zeitlin algebras}   \label{subsection: GZ algebras}
 \subsubsection{Generalities on Gelfan--Zeitlin subalgebras}
 Let us recall the following notion pertaining to a finite multiplicity free tower $(A_k)_{0 \le k \le n}$ of split semisimple algebras over a field $F$. The terminology is from Vershik and Okounkov ~\cite{okounkov-vershik-selecta, okounkov-vershik-JMS-NY}. 
 
 \begin{definition}  The {\em Gelfand--Zeitlin subalgebra} $G_n$  of $A_n$ is the subalgebra generated by the centers of $A_0, A_1, \dots, A_n$.  
 \end{definition}
 
The  Gelfand--Zeitlin subalgebra is a maximal abelian subalgebra of $A_n$ and contains a remarkable family of idempotents indexed by paths on the branching diagram $\B$ of $(A_k)_{0 \le k \le n}$.  
For each $j$ let $\{z_\la : \la \in \La_j\}$ denote the set of minimal central idempotents in $A_j$.
For $k \le n$ and $\mathfrak t$ a path on $\B$ of length $k$, let $F_\mathfrak t = \prod_j z_{\mathfrak t(j)}$.   Then the elements $F_\mathfrak t$ for $\mathfrak t$ of length $k$ are mutually orthogonal minimal idempotents whose sum is the identity;  moreover the sum
of those $F_\mathfrak t$ such that $\mathfrak t(k) = \la$ is $z_\la$.  If $s$ is a path of length $k$ and $t$ is a path of length $\ell$, with $k \le \ell$, then $F_\mathfrak s F_\mathfrak t =  \delta_{\mathfrak s, \mathfrak t[0, k]} F_\mathfrak t$.   Evidently, the set of $F_\mathfrak t$ as $\mathfrak t$ varies over paths of length $k \le n$ generate $G_n$.   Let us call the set $\{F_\mathfrak t\}$ the family of Gelfand--Zeitlin idempotents for $(A_k)_{0 \le k \le n}$.  The properties listed above characterize this family of idempotents:

\begin{lemma} \label{lemma:  characterization of GZ idempotents}
Consider a finite multiplicity free tower $(A_k)_{0 \le k \le n}$ of split semisimple algebras over a field $F$. Let $F'_\mathfrak t$ be  a family of idempotents indexed by paths of length $k \le n$ on the branching diagram $\B$ of $(A_k)_{0 \le k \le n}$ with the following properties:
\begin{enumerate}
\item  For $t$ of length $k$,  $F'_\mathfrak t$ is a minimal idempotent in $A_k$.   The sum of those $F'_\mathfrak t$ such that $t$ has length $k$ and $t(k) = \la$ is $z_\la$.  
\item  If $\mathfrak s$ is a path of length $k$ and $\mathfrak t$ is a path of length $\ell$, with $k \le \ell$, then $F'_s F'_\mathfrak t =  \delta_{\mathfrak s, \mathfrak t[0, k]} F'_\mathfrak t$. 
\end{enumerate}
Then $F'_\mathfrak t = F_\mathfrak t$ for all paths $t$. 
\end{lemma}

\begin{proof} Let  $\mathfrak t$ be a path of length $k \ge 1$ let 
$\mathfrak t' =\mathfrak t[0, k-1]$ and $\la = \mathfrak t(k)$.  It follows from the assumptions that $F'_\mathfrak t = F'_{\mathfrak t'} z_\la$.   Using this, the conclusion $F'_\mathfrak t = F_\mathfrak t$ follows by induction on the length of the path.  
\end{proof}

 \subsubsection{JM elements and GZ subalgebras}
 We return to  our assumptions \ref{assumption: assumptions on coherent tower for section 3}. Suppose that $(L_n)_{n\ge 0}$ is a multiplicative or additive JM family in $(A_n)_{n \ge 0}$.  
 According to Propositions \ref{proposition: triangularity property of JM elements} and \ref{proposition: triangularity property of  additive JM elements},  for each $n \ge 0$,  $\{L_1, \dots, L_n\}$  is a JM family for $A_n$ in the sense of Mathas, with respect to the reverse lexicographic order and any path basis.
 Suppose now, in addition, that Mathas' separation property is satisfied, namely that
 for each $n$,  $\mathfrak t \mapsto (\kappa(\mathfrak t, j))_{1 \le j \le n}$ is injective on $\mathcal T(n)$.
 
 \begin{proposition} \label{proposition: JM elements generate GZ algebra}
Suppose that  for each $k$, $\mathfrak t \mapsto (\kappa(\mathfrak t, j))_{1 \le j \le k}$ is injective on $\mathcal T(k)$.   Then for each $n$, $\{L_1, \dots, L_n\}$  generates the Gelfand--Zeitlin subalgebra of the finite tower $(A_k^F)_{0 \le k \le n}$. 
\end{proposition}  

\begin{proof}  
Fix $n$.  
For $j \le k \le n$, let  $K(j) = \{\kappa(\mathfrak t, j) : \mathfrak t \in \mathcal T(k)\}$; note that
$K(j)$ does not depend on $k$ as long as $j \le k$.  
For $\mathfrak t$ a path on $\B$ of length $k$, define
\newcommand\doubledecker[2]{\genfrac {} {} {0 pt} {} {#1}{#2}}
$$
F'_\mathfrak t = \prod_{j = 1}^k  \prod_{\doubledecker{c \in K(j)}{ c \ne \kappa(\mathfrak t, j)} }\frac{L_j - c}{\kappa(\mathfrak t, j) - c}.
$$
    Then Mathas  ~\cite{mathas-seminormal} shows that  $F'_\mathfrak t$ is a minimal idempotent in $A_k^F$ and the sum of those
    $F'_\mathfrak t$ such that $\mathfrak t$ has length $k$ and $\mathfrak t(k) =  \la$ is $z_\la$.  Moreover, for $j \le k$,  $L_j  F'_\mathfrak t = \kappa(\mathfrak t, j) F'_\mathfrak t$.    It follows from this that if 
    $\mathfrak s$ is a path of length $k$ and $\mathfrak t$ is a path of length $\ell$, with $k \le \ell$, then $F'_s F'_\mathfrak t =  \delta_{\mathfrak s, \mathfrak t[0, k]} F'_\mathfrak t$. 
   Hence, by Lemma \ref{lemma:  characterization of GZ idempotents},  Mathas' idempotents $F'_\mathfrak t$ are the Gelfand--Zeitlin idempotents for the finite tower  $(A_k^F)_{0 \le k \le n}$.  
   This shows that the Gelfand--Zeitlin algebra is contained in the algebra generated by the JM elements; 
   on the other hand, the JM elements are in the linear span of the idempotents $F'_t$, which gives the opposite inclusion.  
\end{proof}

 \section{Framework axioms and a theorem on cellularity} \label{subsection: framework axioms}  We describe the framework axioms and main theorem of ~\cite{GG1}.
Let $R$ be an integral domain with field of fractions $F$.  We consider two towers of $R$--algebras
$$
A_0 \subseteq A_1 \subseteq A_2 \subseteq \cdots, \quad\text{and} \quad Q_0 \subseteq Q_1 \subseteq Q_2 \subseteq \cdots.
$$
The framework axioms of ~\cite{GG1} are the following:

\begin{enumerate}
\item \label{axiom Hn coherent}  $(Q_n)_{n \ge 0}$ is a coherent tower of cellular algebras.
\item \label{axiom: involution on An}  There is an algebra involution $i$  on $\cup_n A_n$ such that $i(A_n) = A_n$.
\item  \label{axiom: A0 and A1}  $A_0 = Q_0 = R$, and $A_1 = Q_1$  (as algebras with involution).
\item  \label{axiom: semisimplicity}
For all $n$,  $A_n^F : = A_n \otimes_R F$   is split semisimple.  
\item \label{axiom:  idempotent and Hn as quotient of An}
 For $n \ge 2$,  $A_n$ contains an  essential idempotent $e_{n-1}$ such that $i(e_{n-1}) = e_{n-1}$ and
\break $A_n/(A_n e_{n-1} A_n) \cong Q_n$,  as algebras with involution.

\item \label{axiom: en An en} For $n \ge 1$,   $e_{n}$ commutes with $A_{n-1}$ and $e_{n} A_{n} e_{n} \subseteq  A_{n-1} e_{n}$.
\item  \label{axiom:  An en}
For $n \ge 1$,  $A_{n+1} 	e_{n} = A_{n} e_{n}$,  and the map $x \mapsto x e_{n}$ is injective from
$A_{n}$ to $A_{n} e_{n}$.
\item \label{axiom: e(n-1) in An en An} For $n \ge 2$,   $e_{n-1} \in A_{n+1} e_n A_{n+1}$.
\end{enumerate}

Say that the pair of towers of algebras  $(Q_k)_{k\ge 0}$ and
$(A_k)_{k\ge 0}$  satisfy the {\em strong framework axioms},  if they satisfy the axioms with
(1) replaced by 

\noindent
\quad\ \ \  (1$'$) \ \    $(Q_n)_{n \ge 0}$ is a strongly coherent tower of cellular algebras.

In the following theorem, point (4) we use the notion of a branching diagram obtained by reflections from another branching diagram.  We refer the reader to ~\cite{GG1},  Section
2.5 for this notion.

\begin{theorem}[\cite{GG1}, Theorem 3.2] \label{main theorem}
 Let $R$ be an integral domain with field of fractions $F$.  Let $(Q_n)_{n\ge 0}$ and
$(A_n)_{n\ge 0}$  be two towers of $R$--algebras satisfying the framework axioms (resp. the strong framework axioms).  Then
\begin{enumerate}
\item $(A_n)_{n\ge 0}$ is a coherent tower of cellular algebras (resp. a strongly coherent tower of cellular algebras). 
\item  For all $n$, the  partially ordered set in the cell datum for $A_n$ can be realized as
$$
\Lambda_n =  \coprod_{\substack{i\leq n\\n-i\text{ even}}} \Lambda_i\spp 0 \times \{n\},
$$
with the following partial order:  Let $\lambda \in \Lambda_i\spp 0$ and $\mu \in \Lambda_j\spp 0$,  with
$i$, $j$, and $n$ all of the same parity.  Then
 $(\lambda, n) > (\mu, n)$ if, and only if,   $i < j$,  or $i = j$ and $\lambda > \mu$ in $\Lambda_i\spp 0$.
 \item  Suppose $n \ge 2$ and $(\la, n)  \in \La_i\spp 0 \times \{n\} \subseteq \La_n$.  Let 
 $\Delta^{(\la, n)}$ be the corresponding cell module.  
 \begin{enumerate}
 \item
  If  $i < n$,  
 then  $\Delta^{(\la, n)} =
  A_{n-1} e_{n-1} \otimes_{A_{n-2}} \Delta^{(\la, n-2)}$.  Moreover,  
 $$(A_n e_{k-1} A_n \ \Delta^{(\la, k)})\otimes_R F = \Delta^{(\la, k)}\otimes_R F.$$
 \item If $i = n$  then   $\Delta^{(\la, n)}$ is a \, $Q_n$ module, and   $A_n e_{n-1} A_n  \ \Delta^{(\la, n)}  = 0$.
 \end{enumerate}
\item The branching diagram $\B$ for $(A_n^F)_{k \ge 0}$ is that obtained by reflections from the branching diagram
$\B_0$ for $(Q_n^F)_{n \ge 0}$.

\end{enumerate}
\end{theorem}

\begin{proof}  The theorem for coherent towers is proved in ~\cite{GG1}. The modification for strongly coherent towers is straightforward. 
\end{proof}

\begin{remark}  At first sight, it may seem that to apply Theorem \ref{main theorem} requires verifying a formidable list of axioms, but in fact 
the theorem  is always easy to apply.  All of the axioms except  (1$'$) and (4) are elementary.   Axiom  (1$'$) is generally a substantial theorem, which however is already available in the literature in many interesting examples.  Axiom (4) can generally be verified by use of Wenzl's method,   applying the Jones basic construction.  For examples, see section \ref{section: examples} of this paper and ~\cite{GG1}, section 5.
\end{remark}
 
   \section{JM elements in algebras arising from the basic construction} \label{section: JM elements}
   
   \begin{theorem} \label{theorem:  JM elements in basic construction algebras 1}
    Consider two towers of $R$--algebras  $(A_n)_{n \ge 0}$ and $(Q_n)_{n \ge 0}$ satisfying the strong  framework axioms of  Section \ref{subsection: framework axioms}.
   Suppose that $\{L_j\spp 0 : j \ge 1\}$  is a multiplicative JM--family for the tower $(Q_n)_{n \ge 0}$, in the sense of Section \ref{subsection: preliminaries: JM}, and that $\{L_n : n \ge 1\}$
   is a family of elements in $(A_n)_{n \ge 0}$ satisfying the following conditions:
   \begin{enumerate}
   \item  $L_n \in A_n$,  and $L_n$  commutes with $A_{n-1}$.  
   \item  $\pi_j(L_j) = L_j\spp 0$, where $\pi_j : A_j \to Q_j$  is the quotient map.
     \item  For each $j \ge 1$, there exists  $\gamma_j \in R^\times$ such that   $$L_j L_{j+1} e_j =  e_j L_j L_{j+1} =\gamma_j  e_j.$$
   \end{enumerate}
     Then $\{L_j : j \ge 1\}$ is a multiplicative JM--family for the tower $(A_n)_{n \ge 0}$. 
  \end{theorem}
  
  \begin{proof}    Write $\La_n\spp 0$ for the partially ordered set  in the cell datum for $Q_n$
  and $\La_n$ for that  in the cell datum for $A_n$.   Recall that $\La_n$ is realized as the set
  of ordered pairs $(\la, n)$,  where $\la \in \La_k\spp 0$ for some $k \le n$  with $n - k$ even.
  For $n \ge 1$  and $\la \in \La_n\spp 0$, let 
 $\alpha(\la) \in R^\times$  be  such that the product $L_1\spp 0 \cdots L_n \spp 0$  acts by the scalar $\alpha(\la)$ on the cell module $\Delta^\la$ of $Q_n$. 
  
To show that $\{L_j : j \ge 1\}$ is a multiplicative JM--family for the tower $(A_n)_{n \ge 0}$, 
  we need only verify point (2) of definition \ref{definition: multiplicative JM family}.   We do this by induction on $n$.  For $n =0$,  we interpret $L_1 \cdots L_n$ to be
  the identity, and we observe that the statement is trivial.  For $n = 1$,  $A_1 = Q_1$, so again there is nothing to prove.   
  Suppose that $n >1$, and that for all $m < n$ and   all  $(\mu, m) \in \La_{m} $,  with $\mu \in \La_k\spp 0$, 
$(L_1 \cdots L_{m})$  acts as the scalar $$\beta((\mu, m)) : =
\gamma_{m-1} \gamma_{m-3} \cdots \gamma_{k+1} 
\alpha(\mu)$$ on the cell module $\Delta^{(\mu, m)} $ of  $A_{m}$.  

  If $\la \in \La_n\spp 0$, then the cell module $\Delta^{(\la, n)}$ is actually the $Q_n$--module $\Delta^\la$, so
  $$
  (L_1 \cdots L_n) y =  (L_1\spp 0 \cdots L_n \spp 0) y = \alpha(\la) y,
  $$
  for $y \in \Delta^{(\la, n)}$.  
  
  Let  $\la \in \La_k\spp 0$ for some $k < n$.    
  Then $\Delta^{(\la, n)} =
  A_{n-1} e_{n-1} \otimes_{A_{n-2}} \Delta^{(\la, n-2)}$.  For $x \in A_{n-1}$ and $y \in 
   \Delta^{(\la, n-2)}$, we have
  $$
  \begin{aligned}
 (L_1 \cdots L_n) x e_{n-1} \otimes y &=  (L_1 \dots L_{n-1}) x L_n e_{n-1} \otimes y \\
 &= x  (L_1 \dots L_{n-1}) L_n e_{n-1} \otimes y \\
 &=    x (L_{n-1} L_n) e_{n-1} \otimes (L_1 \cdots L_{n-2}) y \\
 &= \gamma_{n-1} x e_{n-1} \otimes  \gamma_{n-3} \cdots \gamma_{k+1} \alpha(\la) y\\
 &= \gamma_{n-1} \cdots \gamma_{k+1} \alpha(\la) x e_{n-1} \otimes y,
 \end{aligned}
  $$
   where the first equality is valid since $L_n$ commutes with $A_{n-1}$,   the second follows from  the induction hypothesis and  Proposition \ref{proposition: triangularity property of JM elements} (2),    the third follows  because $L_1 \dots L_{n-2}$  is an element of $A_{n-2}$,  and so commutes with $e_{n-1}$,   and the fourth comes from the induction hypothesis  and hypothesis (3)  of the theorem statement. 
 \end{proof}
  
  \begin{corollary}  If  $\gamma_j$ is independent of $j$,  say $\gamma_j = \gamma$ for all $j$,  then
  $\beta((\la, n)) = \gamma^{(n-k)/2}  \alpha(\la)$ when $\la \in \La_k\spp 0$. 

  \end{corollary}
  
   The additive version of the theorem is the following. The proof is similar.  
   
    \begin{theorem}  \label{theorem:  JM elements in basic construction algebras 2}
    Consider two towers of $R$--algebras  $(A_n)_{n \ge 0}$ and $(Q_n)_{n \ge 0}$ satisfying the strong  framework axioms of  Section \ref{subsection: framework axioms}.
   Suppose that $\{L_j\spp 0 : j \ge 1\}$  is an additive JM--family for the tower $(Q_n)_{n \ge 0}$, in the sense of Section \ref{subsection: preliminaries: JM}, and that $\{L_n : n \ge 1\}$
   is a family of elements in $(A_n)_{n \ge 0}$ satisfying the following conditions:
   \begin{enumerate}
   \item  $L_n \in A_n$,  and $L_n$  commutes with $A_{n-1}$.  
   \item  $\pi_j(L_j) = L_j\spp 0$, where $\pi_j : A_j \to Q_j$  is the quotient map.
    \item  For each $j \ge 1$, there exists  $\gamma_j \in R$ such that  $$(L_j + L_{j+1}) e_j =  e_j ( L_j + L_{j+1}) = \gamma_j e_j.$$ 
       \end{enumerate}
       Then $\{L_j : j \ge 1\}$ is an additive JM--family for the tower $(A_n)_{n \ge 0}$.  \end{theorem}
   
   The additive analogue of the formula for $\beta$ developed in the proof of Theorem 
   \ref{theorem:  JM elements in basic construction algebras 1} is the following. 
    For $n \ge 1$  and $\la \in \La_n\spp 0$, let 
 $d(\la) \in R$  be  such that $L_1\spp 0 + \cdots + L_n \spp 0$  acts by the scalar $d(\la)$ on the cell module $\Delta^\la$ of $Q_n$.  Then for $(\la, n) \in \La_n$, 
 with $\la \in \La_k\spp 0$, 
  $L_1 + \cdots + L_n$
 acts by the scalar
 $$
 \beta((\la, n)) =\gamma_{n-1} + \cdots + \gamma_{k+1}  + d(\la).
 $$
 If $\gamma_j$ is independent of $j$,  say $\gamma_j = \gamma$ for all $j$, then
 $$
\beta((\la, n)) = \frac{n-k}{2} \gamma + d(\la).
$$

    \section{Examples} \label{section: examples}
    
    \subsection{Preliminaries on tangle diagrams}  \label{subsection:  preliminaries on tangle diagrams}
Several of our examples involve {\em tangle diagrams} in the rectangle $\mathcal R = [0, 1] \times [0, 1]$.
Fix points $a_i \in [0, 1]$,  $i \ge 1$,  with $0 < a_1 < a_2 < \cdots$. Write
$\p i = (a_i, 1)$ and $\overline{ \p i} = (a_i, 0)$.

Recall that a {\em knot diagram} means a collection of piecewise smooth closed curves in the plane
which may have intersections and self-intersections, but only simple
transverse intersections.  At each intersection or crossing, one of the
two strands (curves) which intersect is indicated as crossing
over the other.  

An {\em $(n,n)$--tangle diagram}  is a piece of a
knot diagram in $\mathcal R$  consisting of exactly $n$ topological intervals and possibly some number of closed curves, such that:  (1)   the endpoints of the intervals are the points $\p 1, \dots, \p n, \pbar 1, \dots, \pbar n$, and these are the only points of intersection of the family of curves with the boundary of the rectangle, and (2)  each interval intersects the boundary of the rectangle transversally.

An  {\em $(n,n)$--Brauer diagram} is a ``tangle" diagram  containing no closed curves, 
in which information about over and under crossings is ignored.  Two Brauer diagrams are identified if the pairs of boundary points joined by curves is the same in the two diagrams.
By convention, there is a unique $(0, 0)$--Brauer diagram,  the empty diagram with no curves.
For $n \ge 1$,  the number of $(n,n)$--Brauer diagrams is  $(2n-1)!! = (2n-1)(2n-3)\cdots (3)(1)$.

\ignore{
A {\em Temperley--Lieb} diagram is a Brauer diagram without crossings.  For $n \ge 0$,  number of $(n, n)$--Temperley--Lieb diagrams is the Catalan number $\frac{1}{n+1} {2n \choose n}$.
}

For any of these types of diagrams, we call $P = \{\p 1, \dots, \p n, \pbar 1,\dots,  \pbar n\}$ the set of {\em vertices} of the diagram,  $P^+ =    \{\p 1, \dots, \p n\}$ the set of {\em top vertices},  and
$P^- =  \{\pbar 1,\dots,  \pbar n\}$ the set of {\em bottom vertices}.  A curve or {\em strand} in the diagram is called a {\em vertical} or {\em through} strand if it connects a top vertex and a bottom vertex,  and a {\em horizontal} strand if it connects two top vertices or two bottom vertices.

    \subsection{The BMW algebras}  
The BMW algebras were first introduced by Birman and Wenzl ~\cite{Birman-Wenzl} and independently by Murakami ~\cite{Murakami-BMW} as abstract algebras defined by generators and relations.  The version of the presentation given here follows \cite{Morton-Traczyk} and
 ~\cite{Morton-Wassermann}.

\def\hods{\unskip\kern.55em\ignorespaces}
\def\mods{\vskip-\lastskip\vskip4pt}
\def\ods{\vskip-\lastskip\vskip4pt plus2pt}
\def\bods{\vskip-\lastskip\vskip12pt plus2pt minus2pt}

\begin{definition}  \label{definition: BMW algebra}
Let $S$ be a commutative unital ring with invertible elements $\rho$ and $q$ and an element $\delta$ satisfying $\rho\inv - \rho = (q\inv -q)(\delta -1)$.  The {\em Birman--Wenzl--Murakami algebra}
$\bmw n(S; \rho, q, \delta)$ is the unital $S$--algebra 
 with generators $g_i^{\pm 1}$  and
$e_i$ ($1 \le i \le n-1$) and relations:
\begin{enumerate}
\item (Inverses) \hods $g_i g_i\inv = g_i\inv g_i = 1$.
\item (Essential idempotent relation)\hods $e_i^2 = \delta e_i$.
\item (Braid relations) \hods $g_i g_{i+1} g_i = g_{i+1} g_i g_{i+1}$ 
and $g_i g_j = g_j g_i$ if $|i-j|  \ge 2$.
\item (Commutation relations)  \hods $g_i e_j = e_j g_i$  and
$e_i e_j = e_j e_i$  if $|i-j|\ge 2$. 
\item (Tangle relations)\hods $e_i e_{i\pm 1} e_i = e_i$, $g_i
g_{i\pm 1} e_i = e_{i\pm 1} e_i$, and $ e_i  g_{i\pm 1} g_i=   e_ie_{i\pm 1}$.
\item (Kauffman skein relation)\hods  $g_i - g_i\inv = (q - q\inv)(1- e_i )$.
\item (Untwisting relations)\hods $g_i e_i = e_i g_i = \rho\inv e_i$,
and $e_i g_{i \pm 1} e_i = \rho e_i$.
\end{enumerate}
\end{definition}

The BMW algebra $\bmw n$  can also be realized as the algebra of  $(n, n)$--tangle diagrams modulo regular isotopy and the following  {\em Kauffman skein relations:}  
\begin{enumerate}
\item Crossing relation:
$
\quad \inlinegraphic[scale=.6]{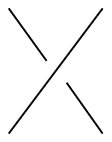} - \inlinegraphic[scale=.3]{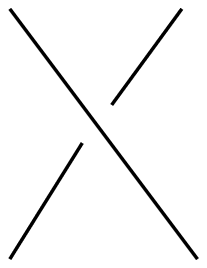} 
\quad = 
\quad
(q\inv - q)\,\left( \inlinegraphic[scale=1.2]{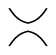} - 
\inlinegraphic[scale=1.2]{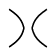}\right).
$
\item Untwisting relation:
$\quad 
\inlinegraphic{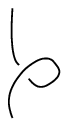} \quad = \quad \rho \quad
\inlinegraphic{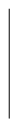} \quad\ \text{and} \quad\ 
\inlinegraphic{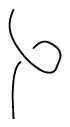} \quad = \quad \rho\inv \quad
\inlinegraphic{vertical_line}. 
$
\item  Free loop relation:  $T\, \cup \, \bigcirc = \delta \, T, $  where $T\, \cup \, \bigcirc$ means the union of a tangle diagram $T$ and a closed loop having no crossings with $T$.
\end{enumerate}
In the tangle picture, 
 $e_j$ and $g_j$ are represented by the following $(n,n)$--tangle diagrams:
$$
e_j =  \inlinegraphic[scale=.7]{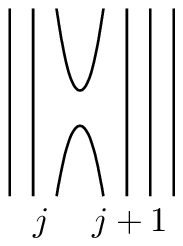}\qquad
g_j =  \inlinegraphic[scale= .7]{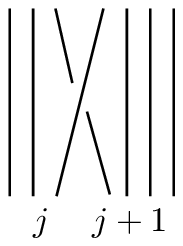} 
$$
The realization of the BMW algebra as an algebra of tangles is   from ~\cite{Morton-Wassermann}.  See ~\cite{GG1}, Section 5.4 for more details. 

  The quotient of 
the BMW algebra $\bmw n(S; \rho, q, \delta)$ by the ideal $J$ generated by $e_{n-1}$ is the Hecke algebra $H_n(S; q^2)$.  If $\pi_n$ denotes the quotient map $\pi_n : \bmw n \to
\bmw n/J$,  take $T_i = \pi_n(q\ g_i)$  to obtain an isomorphism with the Hecke algebra as presented in Example \ref{example: hecke algebras 1}.

 The generic  ground ring for the BMW algebras is
 $$
 R = \Z[\rhobold^{\pm1}, \qbold^{\pm1}, \deltabold]/\langle \rhobold\inv - \rhobold = (\qbold\inv - \qbold)(\deltabold -1) \rangle,
 $$
 where $\rhobold$, $\qbold$, and $\deltabold$ are indeterminants over $\Z$.  
$R$ is an integral domain whose field of fractions is $F \cong \Q(\rhobold, \qbold)$  (with $\deltabold =  
 (\rhobold\inv - \rhobold)/(\qbold\inv - \qbold) + 1$ in $F$.)  
Write $\bmw n$ for 
 $\bmw n(R; \rhobold, \qbold, \deltabold)$  and $H_n$ for $H_n(R; \qbold^2)$.
 It is shown in ~\cite{GG1}, Section 5.4,  that the pair of towers  $(\bmw n)_{n \ge 0}$  and $(H_n)_{n \ge 0}$ satisfy the framework axioms of  Section \ref{subsection: framework axioms}.
 In fact, by Example \ref{example: hecke algebras 1},  the tower of Hecke algebras is
 strongly coherent, so the pair satisfies the strong version of the framework axioms.  Consequently, by Theorem \ref{main theorem}, the sequence of BMW algebras is a strongly coherent tower of cellular algebras.  The partially ordered set $\La_n$  in the cell datum of
 $\bmw n$ is the set of pairs $(\la, n)$,  with $\la$ a Young diagram of size $k \le n$ with
 $n-k$ even.   The set of paths $\mathcal T((\la, n))$ can be identified with up--down tableaux
 of length $n$ and shape $\la$,  see ~\cite{Enyang2}.

 The following analogue of Jucys--Murphy elements for the BMW algebras were introduced 
 by Leduc and Ram ~\cite{leduc-ram}  and Enyang ~\cite{Enyang2}.  
 Define $L_1 = 1$ and $L_{j+1} =  g_j L_j g_j$  for $j \ge 1$.  
 (Thus, for example, $L_5 = g_4 g_3 g_2 g_1^2 g_2 g_3 g_4$.) 
  The involution on $\bmw n$ is the unique algebra involution taking $e_i \mapsto e_i$  and $g_i \mapsto g_i$;  it leaves each
 $L_j$ invariant.  One can check algebraically that $L_n$ commutes with the generators of 
 $\bmw {n-1}$,   but this is far easier to see using  the geometric realization of $\bmw n$. 
 In fact,  in the geometric picture,  $L_n$ is represented by the braid in which the $n$--th strand wraps once around the first  $(n-1)$ strands. 
 
  Let $L_j\spp 0$ denote the classical JM elements in the Hecke algebras $H_n$, as defined in Example  \ref{example: JM elements in Hecke algebra}.  Then we have $\pi_n(L_j) = L_j\spp 0$  for $1 \le j \le n$; this follows because  $\pi_n(L_1) = 1$  and  $\pi_n(L_{j+1}) =  \qbold^{-2} T_j \pi_n(L_j)  T_j$.
 (This is the correct recursion, because the Hecke algebra parameter $q$  has been replaced by $\qbold^2$.)
 One can check, using algebraic relations or by using tangle diagrams, that for all $j \ge 1$, 
 $$
L_j L_{j+1} e_j =  e_j L_j L_{j+1} =
 \rho^{-2} e_j.
 $$
 (The factor of $ \rho^{-2} $  comes from two applications of the untwisting relation (2) above.) 

It now follows from Theorem \ref{theorem:  JM elements in basic construction algebras 1} that
$\{L_j : j \ge 0\}$  is a multiplicative JM--family in $(\bmw n)_{n\ge 0}$, with
$L_1 \dots L_n$  acting by $$ \beta((\la, n)) := \rho^{-(n-k)} \alpha(\la)$$  on the cell module $\Delta^{(\la, n)}$, if
$\la$ is a Young diagram of size $k$.  By Proposition \ref{proposition: triangularity property of JM elements},  the action of the elements $L_j$ on the basis of $\Delta^{(\la, n)}$ labelled by up--down tableaux is triangular:
 \begin{equation} \label{equation: triangular action of Lj bmw case}
   L_j a_\mathfrak t^\la = \kappa(j, \mathfrak t)\  a_\mathfrak t^\la + \sum_{\mathfrak s\, \succ \, \mathfrak t} r_\mathfrak s a_\mathfrak s^\la,
 \end{equation}
 with $\kappa(j, \mathfrak t) = \displaystyle\frac{\beta(\mathfrak t(j))}{\beta(\mathfrak t(j-1))}$,
   for some elements $r_\mathfrak s \in R$, depending on $j$ and $\mathfrak t$. 
   Moreover, if $\mathfrak t(j) = (\nu, j)$  and $\mathfrak t(j-1) = (\mu, j-1)$, then
   $|\nu| = |\mu| \pm 1$.  If   $|\nu| = |\mu| + 1$ and $\nu \setminus \mu = x$,  then
   $$
   \kappa(j, \mathfrak t) = \displaystyle\frac{\beta((\nu,j))}{\beta((\mu,j-1))} = 
   \displaystyle\frac{\alpha(\nu)}{\alpha(\mu)} = \qbold^{2 \kappa(x)},
     $$
     where $\kappa(x)$ is the content of $x$, namely the column of $x$ minus the row of $x$.  
    If   $|\nu| = |\mu| - 1$ and $\mu \setminus \nu = x$,  then
   $$
   \kappa(j, \mathfrak t) = \displaystyle\frac{\beta((\nu,j))}{\beta((\mu,j-1))} = 
  \rho^{-2} \displaystyle\frac{\alpha(\nu)}{\alpha(\mu)} = \rho^{-2}\qbold^{- 2 \kappa(x)}.
     $$
This recovers Theorem 7.8 of Enyang ~\cite{Enyang2}.\footnote{The theorem is stated in 
 ~\cite{Enyang2} with dominance order rather than lexicographic order, but it appears that the proof only yields the statement with lexicographic order.} 
 
 \subsection{The Brauer algebras}
 The Brauer algebras were introduced by  Brauer~\cite{Brauer} as a device
for studying the invariant theory of orthogonal and symplectic groups.

 Let $S$ be a commutative ring with identity, with a distinguished element $\delta$.
The Brauer algebra $\br_n(S, \delta)$ is the free $S$--module with basis the set of $(n, n)$--Brauer diagrams,  with multiplication defined as follows.
The product of two Brauer diagrams is defined
to be a certain multiple of another Brauer diagram.  Namely, given two
Brauer diagrams $a, b$, first ``stack" $b$ over $a$; the result is a planar tangle that may contain some number of closed curves.    Let $r$ denote the number of closed curves, and let $c$ be the Brauer
diagram obtained by removing all the closed curves.  Then
$
a b = \delta^r c.
$

\begin{definition}\rm
For $n \ge 1$, the {\em Brauer algebra} $\br_n(S, \delta)$ over $S$ with parameter $\delta$ is the free $S$-module with basis the set of                     
$(n,n)$-Brauer diagrams, with the bilinear product determined by the
multiplication of Brauer diagrams.  In particular, $\br_0(S, \delta) = S$.
\end{definition}

Note that the Brauer diagrams with only vertical strands are in
bijection with permutations of $\{1, \dots, n\}$, and that the
multiplication of two such diagrams coincides with the multiplication of
permutations.  Thus the  Brauer algebra contains the group algebra $S\S_n$ of
the permutation group $\mathfrak S_n$.  The identity element of the Brauer algebra is the diagram corresponding to the trivial permutation.  We will note below that $S \S_n$ is also a quotient of $\br_n(S, \delta)$.  

The involution $i$ on $(n, n)$--Brauer diagrams which reflects a diagram in the axis $y = 1/2$
extends linearly to an algebra involution of $\br_n(S, \delta)$.

Let $e_j$ and $s_j$ denote the $(n, n)$--Brauer diagrams:
$$
e_j =  \inlinegraphic[scale=.7]{ordinary_E_j}\qquad
s_j =  \inlinegraphic[scale= .7]{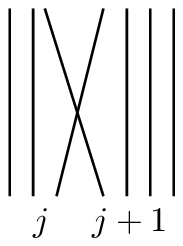} 
$$
Note that $e_j^2 = \delta e_j$,  so $e_j$ is an essential idempotent if $\delta \ne 0$, and nilpotent if $\delta = 0$.
We have $i(e_j) = e_j$ and $i(s_j) = s_j$.  
It is easy to see that $e_1, \dots, e_{n-1}$ and $s_1, \dots, s_{n-1}$ generate $\br_n(S, \delta)$ as an algebra.

The products $a b$ and $b a$  of two Brauer diagrams have  at most as many through strands as $a$.  Consequently, the span of diagrams with fewer than $n$ through strands is an ideal 
$J$ in $\br_n(S, \delta)$.  
 The
ideal  $J $  is generated by $e_{n-1}$.  We have   $\br_n(S, \delta)/J \cong S\S_n$, as algebras with involutions.

The generic ground ring for the Brauer algebras is  $R = \Z[\deltabold]$,  where $\deltabold$ is an indeterminant.    Let $F = \Q(\deltabold)$ denote the field of fractions of $R$.  Write $\br_n = \br_n(R, \deltabold)$.

 It is shown in ~\cite{GG1}, Section 5.2,  that the pair of towers  $(\br_n)_{n \ge 0}$  and $(R\S_n)_{n \ge 0}$ satisfy the framework axioms of  Section \ref{subsection: framework axioms}.
 In fact, since the symmetric group algebra is a specialization of the Hecke algebra, 
  the tower of symmetric group algebras is
 strongly coherent, so the pair satisfies the strong version of the framework axioms.  Consequently, by Theorem \ref{main theorem}, the sequence of Brauer algebras is a strongly coherent tower of cellular algebras.  As for the BMW algebras, the partially ordered set $\La_n$  in the cell datum of
 $\br_n$ is the set of pairs $(\la, n)$,  with $\la$ a Young diagram of size $k \le n$ with
 $n-k$ even.   The set of paths $\mathcal T((\la, n))$ can be identified with up--down tableaux
 of length $n$ and shape $\la$.
 
  We need to recall the Jucys--Murphy elements for the symmetric group algebras, which can be defined inductively by $L_1\spp 0 = 0$,  $L_{j+1}\spp 0 =  s_j L_{j} s_j + s_j$.  Thus, for example,  $L_5\spp 0 = (1, 5) + (2, 5) + (3, 5) + (4, 5)$.     One has $L_j \spp 0 \in R \S_j$, and $L_j\spp 0$ commutes with $R\S_{j-1}$.  $L_1 \spp 0 + \cdots + L_n \spp 0$ is central
  in $R \S_n$  and acts as the scalar $\alpha(\la) = \sum_{x \in \la} \kappa(x)$ on the cell module
  $\Delta^\la$.  Here, $\la$ is  a Young diagram of size $n$ and for a cell $x$ of $\la$,
  $\kappa(x)$ is the content of $x$,  namely the column co-ordinate   minus the row co-ordinate of $x$. In particular $\{L_j\spp 0 : j \ge 0\}$  is an additive JM--family in the sense of Definition
  \ref{definition: additive JM family}.
  
  The following analogues of Jucys-Murphy elements for the Brauer algebras were introduced by Nazarov ~\cite{Nazarov}.  Let $L_1 = 0$ and $L_{j+1} = s_j L_j s_j +  s_j - e_j  $.  Observe that $\pi_n(L_j) = L_j\spp 0$  for $1 \le j \le n$, where $\pi_n : \br_n \to R\S_n$ is the quotient map.  
  Evidently,  $L_n \in \br_n$.  
    By
 ~\cite{Nazarov}, Proposition  2.3,  $L_n$ commutes with $\br_{n-1}$, and for all $j \ge 1$, 
 $$(L_j + L_{j+1}) e_j = e_j (L_j + L_{j+1}) = (1- \delta) e_j.$$ 
 
 It now follows from Theorem \ref{theorem:  JM elements in basic construction algebras 2} that
$\{L_j : j \ge 0\}$  is an additive JM--family in $(\br_n)_{n\ge 0}$, with
$L_1 + \dots + L_n$  acting by $$ \beta((\la, n)) := \frac{n-k}{2} (1- \delta) +  \alpha(\la)$$  on the cell module $\Delta^{(\la, n)}$, if
$\la$ is a Young diagram of size $k$. 

 By Proposition \ref{proposition: triangularity property of  additive JM elements},  the action of the elements $L_j$ on the basis of $\Delta^{(\la, n)}$ labelled by up--down tableaux is triangular:
 \begin{equation} \label{equation: triangular action of Lj bmw case}
   L_j a_\mathfrak t^\la = \kappa(j, \mathfrak t)\  a_\mathfrak t^\la + \sum_{\mathfrak s\, \succ \, \mathfrak t} r_\mathfrak s a_\mathfrak s^\la,
 \end{equation}
  with $\kappa(j, \mathfrak t) = {\beta(\mathfrak t(j))}- {\beta(\mathfrak t(j-1))}$,
   for some elements $r_\mathfrak s \in R$, depending on $j$ and $\mathfrak t$. 
    Moreover, if $\mathfrak t(j) = (\nu, j)$  and $\mathfrak t(j-1) = (\mu, j-1)$, then
   $|\nu| = |\mu| \pm 1$.  If   $|\nu| = |\mu| + 1$ and $\nu \setminus \mu = x$,  then
$$
   \kappa(j, \mathfrak t) = {\beta((\nu,j))}-{\beta((\mu,j-1))} = 
{\alpha(\nu)}-{\alpha(\mu)} =  \kappa(x).
     $$
  If   $|\nu| = |\mu| - 1$ and $\mu \setminus \nu = x$,  then
   $$
   \kappa(j, \mathfrak t) = {\beta((\nu,j))}-{\beta((\mu,j-1))} = 
(1-\delta) +  {\alpha(\nu)}- {\alpha(\mu)} = (1-\delta)  -\kappa(x) .
     $$
This recovers Theorem 10.7 of Enyang ~\cite{Enyang2}.\footnote{The same caution about lexicographic order versus dominance order applies here, as in the BMW case.}

\subsection{Cyclotomic BMW algebras}
The {\em cyclotomic Birman--Wenzl--Murakami  algebras} are BMW analogues of cyclotomic Hecke algebras  ~\cite{ariki-koike, ariki-book}.  The cyclotomic BMW algebras  were defined by 
H\"aring--Oldenburg in ~\cite{H-O2}  and have recently been studied by three groups of mathematicians:
Goodman and   Hauschild--Mosley ~\cite{GH1, GH2, GH3,  goodman-2008,  goodman-admissibility},  Rui, Xu, and Si  ~\cite{rui-2008, rui-2008b},   and Wilcox and Yu  ~\cite{Wilcox-Yu, Wilcox-Yu2, Wilcox-Yu3, Yu-thesis}.  

\subsubsection{Definition of cyclotomic BMW algebras}
\begin{definition} Fix an integer $r \ge 1$.   A {\em ground ring} 
$S$ is a commutative unital ring with parameters $\rho$, $q$, $\delta_j$   ($j \ge 0$),  and
$u_1, \dots, u_r$, with    $\rho$, $q$,   and $u_1, \dots, u_r$  invertible,  and with $\rho\inv - \rho=   (q\inv -q) (\delta_0 - 1)$.
\end{definition}

\begin{definition} \label{definition:  cyclotomic BMW}
Let $S$ be a ground ring with
parameters $\rho$, $q$, $\delta_j$   ($j \ge 0$),  and
$u_1, \dots, u_r$.
The {\em cyclotomic BMW algebra}  $\bmw{n, S, r}(u_1, \dots, u_r)$  is the unital $S$--algebra
with generators $y_1^{\pm 1}$, $g_i^{\pm 1}$  and
$e_i$ ($1 \le i \le n-1$) and relations:
\begin{enumerate}
\item (Inverses)\hods $g_i g_i\inv = g_i\inv g_i = 1$ and 
$y_1 y_1\inv = y_1\inv y_1= 1$.
\item (Idempotent relation)\hods $e_i^2 = \delta_0 e_i$.
\item (Affine braid relations) 
\begin{enumerate}
\item[\rm(a)] $g_i g_{i+1} g_i = g_{i+1} g_ig_{i+1}$ and 
$g_i g_j = g_j g_i$ if $|i-j|  \ge 2$.
\item[\rm(b)] $y_1 g_1 y_1 g_1 = g_1 y_1 g_1 y_1$ and $y_1 g_j =
g_j y_1 $ if $j \ge 2$.
\end{enumerate}
\item[\rm(4)] (Commutation relations) 
\begin{enumerate}
\item[\rm(a)] $g_i e_j = e_j g_i$  and
$e_i e_j = e_j e_i$  if $|i-
j|
\ge 2$. 
\item[\rm(b)] $y_1 e_j = e_j y_1$ if $j \ge 2$.
\end{enumerate}
\item[\rm(5)] (Affine tangle relations)\vadjust{\vskip-2pt\vskip0pt}
\begin{enumerate}
\item[\rm(a)] $e_i e_{i\pm 1} e_i = e_i$,
\item[\rm(b)] $g_i g_{i\pm 1} e_i = e_{i\pm 1} e_i$ and
$ e_i  g_{i\pm 1} g_i=   e_ie_{i\pm 1}$.
\item[\rm(c)\hskip1.2pt] For $j \ge 1$, $e_1 y_1^{ j} e_1 = \delta_j e_1$. 
\vadjust{\vskip-
2pt\vskip0pt}
\end{enumerate}
\item[\rm(6)] (Kauffman skein relation)\hods  $g_i - g_i\inv = (q - q\inv)(1- e_i)$.
\item[\rm(7)] (Untwisting relations)\hods $g_i e_i = e_i g_i = \rho \inv e_i$
 and $e_i g_{i \pm 1} e_i = \rho  e_i$.
\item[\rm(8)] (Unwrapping relation)\hods $e_1 y_1 g_1 y_1 = \rho e_1 = y_1 
g_1 y_1 e_1$.
\item[\rm(9)](Cyclotomic relation) \hods $(y_1 - u_1)(y_1 - u_2) \cdots (y_1 - u_r) = 0$.
\end{enumerate}
\end{definition}

Thus,  a cyclotomic BMW algebra is the quotient of the affine BMW algebra  ~\cite{GH1}, by the cyclotomic relation $(y_1 - u_1)(y_1 - u_2) \cdots (y_1 - u_r) = 0$. 

The cyclotomic BMW algebra has a unique algebra involution $i$ fixing each of the generators.

\subsubsection{Geometric realization}
 It is shown in ~\cite{GH3}  and in ~\cite{Wilcox-Yu3}  that
 the cyclotomic BMW algebra has a geometric realization as the  ``cyclotomic Kauffman tangle (KT)  algebra,"  assuming admissibility conditions on the ground ring (see below).
The cyclotomic KT algebra is described in terms of  ``affine tangle diagrams,"  which are  just  ordinary  tangle diagrams with a distinguished
 vertical strand connecting $\p 1$  and $\pbar  1$, as in the following figure.
 $$
\inlinegraphic[scale=1.5]{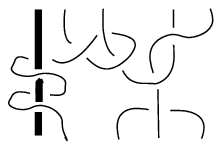}
$$
The { cyclotomic KT algebra} is the algebra of affine tangle diagrams, modulo regular isotopy, Kauffman skein relations,  and a 
 cyclotomic skein relation,  which is a ``local" version of the cyclotomic relation of 
 Definition \ref{definition:  cyclotomic BMW} (9).   See
 ~\cite{GH2}  for the precise definition.

In the geometric realization,  the generators $g_i$,  $e_i$,  and $x_1 = \rho\inv y_1$  are
represented  by the following affine tangle diagrams:
$$
x_1 = \inlinegraphic[scale= .7]{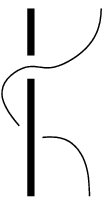}
\qquad
g_i =  \inlinegraphic[scale=.6]{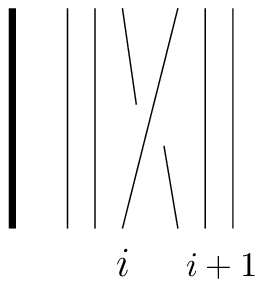}\qquad
e_i =  \inlinegraphic[scale= .7]{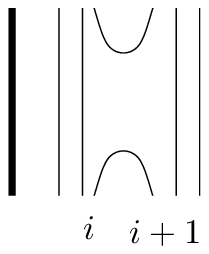}.
$$
In the geometric picture,  the algebra involution $i$ is given on the level of affine tangle diagrams by the map that flips an affine tangle diagram over the horizontal line $y = 1/2$. 

 \subsubsection{Admissibility}
 The cyclotomic BMW algebras can be defined over an arbitrary  ground ring.   
 However, unless the parameters satisfy certain restrictions,  the element $e_1$ is forced to be zero and the algebras collapse to a specialization of the cyclotomic Hecke algebras.   In order to understand the algebras and the restrictions on the parameters, it crucial first to focus on the following ``optimal" situation:  
 
\begin{definition}    Let $S$ be an integral  ground ring with
parameters $\rho$, $q$, $\delta_j$ ($j \ge 0$)  and $u_1, \dots, u_r$,  with $q - q\inv \ne 0$.
One says that $S$ is {\em admissible} (or that the parameters are {\em admissible})   if  \break $\{e_1, y_1 e_1, \dots, y_1^{r-1} e_1\} \subseteq  \bmw{2, S, r}$ is linearly independent over $S$.
\end{definition}

We are also going to restrict our attention to the case that the ground ring is an integral domain, and $q-q\inv \ne 0$.

It is shown by Wilcox and Yu in ~\cite{Wilcox-Yu} that admissibility is equivalent to finitely many (explicit) polynomial conditions on the parameters.  Moreover,  these relations give $\rho$  and $(q - q\inv) \delta_j$  as  Laurent polynomials
in the remaining parameters $q, u_1, \dots, u_r$;    see ~\cite{Wilcox-Yu}  and ~\cite{GH3}  for details.
An alternative set of explicit conditions on the parameters was proposed by Rui and Xu
~\cite{rui-2008}.  It has been shown in ~\cite{goodman-admissibility}  that the conditions of Rui and Xu are also equivalent 
to admissibility
(assuming the ground ring is integral and $q-q\inv \ne 0$). 

Finally, it has been  shown that  the structure of cyclotomic BMW algebras with non--admissible parameters can be derived from the admissible case ~\cite{semi-admissible}.

\subsubsection{Generic ground ring}
There is a universal admissible integral ground ring $R$  for cyclotomic BMW algebras, which is a little more complicated to describe than the generic ground rings for the other algebras we have encountered.  We refer to ~\cite{GH3}, Theorem 3.19 for details.  Suffice it to say that the field of fractions $F$ of $R$ is    $\Q(\qbold, \ubold_1, \dots, \ubold_r)$,
where $\qbold, \ubold_1, \dots, \ubold_r$ are algebraically independent indeterminants over $\Q$;   the remaining parameters are given by certain Laurent polynomials in $\qbold$, $ \ubold_1, \dots, \ubold_r$,
and $(\qbold - \qbold\inv)\inv$, and $R$ is the subring of $F$ generated by all the parameters.  
Any other admissible integral ground ring $S$ is a module over $R$,  and
$\bmw {n, S, r} \cong \bmw {n, R, r} \otimes_R S$.   We will write $\bmw n$  for $\bmw {n, R, r}$.   

\subsubsection{Cyclotomic BMW algebras and cyclotomic Hecke algebras}
We recall the definition of the affine and cyclotomic Hecke algebras,   see ~\cite{ariki-book}.

\begin{definition}
Let $S$ be a commutative unital ring with an invertible element $q$.  The {\em affine Hecke algebra} 
$\ahec{n,S}(q)$ 
over $S$
is the $S$--algebra with generators $T_0,  T_1,  \dots, T_{n-1}$, with relations:
\begin{enumerate}
\item  The generators $T_i$ are invertible,  satisfy the braid relations,  and   the Hecke relations $(T_i - q)(T_i + q) = 0$. 

\item  The generator $T_0$ is invertible,  $T_0 T_1 T_0 T_1 = T_1 T_0 T_1  T_0$  and $T_0$ commutes with $T_j$  for $j \ge 2$.
\end{enumerate}
Let $u_1,  \dots, u_r$  be additional  elements in $S$.   The {\em cyclotomic Hecke algebra}
   $\hec{n, S, r}(q; u_1,  \dots, u_r)$  is the quotient of the affine Hecke algebra $\ahec{n, S}(q)$ by the polynomial relation    $(T_0 - u_1) \cdots (T_0 - u_r) = 0$.
\end{definition}

We remark that since the generator $T_0$ can be rescaled by an arbitrary invertible element of $S$, only the ratios of the parameters $u_i$ have invariant significance in the definition of the cyclotomic Hecke algebra. The cyclotomic Hecke algebra has a unique algebra involution $i$ leaving each generator invariant.  By ~\cite{ariki-koike},  the cyclotomic Hecke algebras $\hec{n, S, r}$  are free $S$--modules of rank $r^n n !$ and $\hec{n, S, r}$  imbeds in 
 $\hec{n + 1, S, r}$. 

The cyclotomic Hecke algebras were shown to be cellular algebras in ~\cite{Graham-Lehrer-cellular}. In ~\cite{dipper-james-mathas}, a cellular basis was given generalizing the Murphy basis of the ordinary Hecke algebra.  The partially ordered set $\La_n \spp 0$
in the cell datum for 
$\hec{n, S, r} = \hec{n, S, r}(q; u_1, \dots, u_r)$  is the set of $r$--tuples of Young diagrams with total size $n$, ordered by dominance.   For each $\labold \in \La_n \spp 0$,   the
index set $\mathcal T(\labold)$ in the cell datum is the set of standard tableaux of shape
$\labold$;  this has the usual meaning:  fillings with the numbers $1, \dots,  n$,  so that
the numbers increase in each row and column (separately in each component Young diagram).    The cyclotomic Hecke algebras are generically split semisimple;  in the semisimple case,  the branching diagram has vertices at level $n$ labelled by all $r$--tuples of Young diagrams of total size $n$, and $\labold \nearrow \mubold$  if $\mubold$ is obtained from $\labold$ by adding one box in one component of $\labold$.   Standard tableaux  of shape $\labold$  can be identified with paths on the generic branching diagram
from $\emptyset$  (the $r$--tuple of empty Young diagrams)  to $\labold$.  

By results of Ariki and Mathas (\cite{ariki-mathas}, Proposition 1.9)  and Mathas 
~\cite{mathas-2009}, the sequence of cyclotomic Hecke algebras $(\hec{n, S, r})_{n\ge 0}$
is a strongly coherent tower of cellular algebras. 

Let $J$ be the ideal in $\bmw n = \bmw {n, R, r}$ generated by $e_{n-1}$.  It is not hard to show that the quotient $\bmw n/ J$ is isomorphic to $\hec{n, R, r}(\qbold^2; \ubold_1, \dots, \ubold_r)$.
If $\pi_n$ denotes the quotient map $\pi_n : \bmw n \to
\bmw n/J$,  take $T_j = \pi_n(q\ g_i)$  for $j \ge 1$,  and
$T_0 = \pi_n(y_1)$   to obtain an isomorphism with the cyclotomic Hecke algebra as presented above.  We will write $\hec n$  for $\hec{n, R, r}(\qbold^2; \ubold_1, \dots, \ubold_r)$.

 It is shown in ~\cite{GG1}, Section 5.5, that the pair of towers of algebras
 $(\bmw n)_{n \ge 0}$  and $(\hec n)_{n \ge 0}$  satisfies the framework axioms of Section 
\ref{subsection: framework axioms}.  Since the sequence of Hecke algebras is strongly coherent, the pairs satisfies the strong version of the framework axioms.  Therefore, it follows from Theorem \ref{main theorem} that  the sequence of cyclotomic BMW  algebras is a strongly coherent tower of cellular algebras. 

The partially ordered set $\La_n$  in the cell datum of
 $\bmw n$ is the set of pairs $(\labold, n)$,  with $\labold$ an $r$-tuple of  Young diagrams of total size $k \le n$ with
 $n-k$ even.   The set of paths $\mathcal T((\labold, n))$ can be identified with up--down tableaux
 of length $n$ and shape $\labold$,  that is  sequences of $r$--tuples of Young diagrams in which each successive $r$--tuple is obtained from the previous one by either adding or removing one box from one component Young diagram.  
 
 \subsubsection{JM elements for cyclotomic BMW and Hecke algebras}
 In the cyclotomic Hecke algebra $\hec {n, S, r} = \hec {n, S, r}(q; u_1, \dots, u_r)$,
 define $L_1 \spp 0 = T_0$  and $L_{j+1}\spp 0 =  q\inv T_j L_j\spp 0T_j$  for $j \ge 1$.  
Then  $L_n\spp 0 \in \hec {n, S, r}$,   $L_n\spp 0$ is invariant under the involution on $\hec {n, S, r}$,
and $L_n\spp 0$ commutes with $ \hec {n-1, S, r}$. The product $L_1\spp 0 \cdots L_n \spp0$ is central in
$ \hec {n, S, r}$.  

 For an $r$--tuple of Young diagrams
$\labold$ of total size $n$ and a cell $x \in \labold$, the multiplicative content of the cell is
$$\kappa(x) = u_j q^{b - a}$$  if $x$ is in row $a$ and column $b$ of the $j$--th component of $\labold$.   For a standard tableau $\mathfrak t$  of shape $\labold$,  and $1 \le j \le n$,  let  $\kappa(j, \mathfrak t) =
\kappa(x)$, where $x$ is the cell occupied by $j$ in $\mathfrak t$.   Let $\{ a_{\mathfrak t}^\labold\}$  be the Murphy type basis of the cell module $\Delta^\labold$  indexed by standard tableaux of shape $\labold$.   Then $L_j\spp 0$ acts by
 \begin{equation} \label{equation: triangular action of Lj cyclotomic Hecke case}
   L_j \spp 0a_\mathfrak t^\la = \kappa(j, \mathfrak t)\  a_\mathfrak t^\la + \sum_{\mathfrak s\, \rhd \, \mathfrak t} r_\mathfrak s a_\mathfrak s^\la,
 \end{equation}
 where the sum is over standard tableaux $\mathfrak s$ greater than $\mathfrak t$ in dominance order  (hence in lexicographic order).   These results are from
 ~\cite{jantzen-sum-formula}, Section 3.   It follows that the product
 $L_1\spp 0 \cdots L_n\spp 0$  acts as the scalar $\alpha(\labold) = \prod_{x \in \labold}  \kappa(x)$ on
 the cell module $\Delta^\labold$.  Thus  $\{L_n\spp 0 : n \ge 0\}$ is a multiplicative  JM--family in the strongly coherent tower of cellular algebras  $(\hec {n, S, r})_{n \ge 0}$.
 
 Define elements $L_j$  in the cyclotomic BMW algebras 
 $\bmw n = \bmw {n, R, r}(\qbold; \ubold_1, \dots, \ubold_r)$ over the generic integral admissible ground ring $R$ by $L_1 = y_1$,  $L_{j+1} =  g_j L_j g_j$  for $j \ge 1$.  These are the same as the elements $y_j$  in ~\cite{GH3}.    We have $L_n \in \bmw n$  and
 $L_n$ commutes with $\bmw  {n-1}$.  One can verify that $L_j L_{j+1} e_j = e_j L_j L_{j+1} = e_j$.  The computations can be done at the level of the affine BMW algebras,using the algebraic relations or using affine tangle diagrams.  
 
 We have $\pi_n(L_1) = T_0 = L_1 \spp 0$,  and
 $\pi_n(L_{j+1}) =  \qbold^{-2} T_j  \pi_n(L_j)  T_j$.  Hence  $\pi_n(L_j)$ satisfy the recursion
 for $L_j\spp 0$ in $\hec n = \hec {n, R, r}(\qbold^2; \ubold_1, \dots, \ubold_r)$.  
  
It now follows from Theorem \ref{theorem:  JM elements in basic construction algebras 1} that
$\{L_j : j \ge 0\}$  is a multiplicative JM--family in $(\bmw n)_{n\ge 0}$, with the product
$L_1 \dots L_n$  acting by $$ \beta((\labold, n)) :=  \alpha(\labold)$$  on the cell module $\Delta^{(\labold, n)}$, if
$\labold$ is an $r$--tuple of Young diagrams of  total size $k$.  By Proposition \ref{proposition: triangularity property of JM elements},  the action of the elements $L_j$ on the basis of $\Delta^{(\labold, n)}$ labelled by up--down tableaux is triangular:
 \begin{equation} \label{equation: triangular action of Lj bmw case}
   L_j a_\mathfrak t^\la = \kappa(j, \mathfrak t)\  a_\mathfrak t^{(\labold, n)} + \sum_{\mathfrak s\, \succ \, \mathfrak t} r_\mathfrak s a_\mathfrak s^{(\labold, n)},
 \end{equation}
 with $\kappa(j, \mathfrak t) = \displaystyle\frac{\beta(\mathfrak t(j))}{\beta(\mathfrak t(j-1))}$,
   for some elements $r_\mathfrak s \in R$, depending on $j$ and $\mathfrak t$. 
   Moreover, if $\mathfrak t(j) = (\nubold, j)$  and $\mathfrak t(j-1) = (\mubold, j-1)$, then
   $|\nubold| = |\mubold| \pm 1$.  If   $|\nubold| = |\mubold| + 1$ and $\nubold \setminus \mubold = x$, 
   where $x$ is the cell in row $a$ and column $b$ of the $\ell$--th  component of $\nubold$, 
    then
   $$
   \kappa(j, \mathfrak t)  = 
   \displaystyle\frac{\alpha(\nubold)}{\alpha(\mubold)} =  \kappa(x) = \ubold_\ell \qbold^{2(b -a)}.
     $$
    If   $|\nubold| = |\mubold| - 1$ and $\mubold \setminus \nubold = x$,  then
   $$
   \kappa(j, \mathfrak t) =
\displaystyle\frac{\alpha(\nubold)}{\alpha(\mubold)} =  \kappa(x)\inv = \ubold_\ell\inv \qbold^{-2(b - a)} .
     $$
This recovers Theorem 3.17 of Rui and Si ~\cite{rui-2008b}.  

\subsection{Degenerate cyclotomic BMW algebras (cyclotomic Nazarov Wenzl algebras)}
Degenerate affine BMW algebras were introduced by Nazarov ~\cite{Nazarov}  under the name {\em affine Wenzl algebras}.    
The   cyclotomic  quotients of these algebras were introduced by Ariki, Mathas, and Rui in \cite{ariki-mathas-rui}  and studied further by Rui and Si in ~\cite{rui-si-degnerate}, under the name {\em cyclotomic Nazarov--Wenzl algebras}.
We propose to refer to these algebras as degenerate  affine  (resp. degenerate cyclotomic)  BMW algebras instead,   to bring the terminology in line with that used for degenerate affine and  cyclotomic Hecke algebras.

\subsubsection{Definition of the degenerate cyclotomic BMW algebras}
Fix a positive integer $n$ and a commutative ring $S$ with
multiplicative identity.  Let  $\Omega=\{\omega_a: \ a\ge0\} $  be a sequence of elements of $S$.

\def\W{W}
\def\Waff{\W^{\text{aff}}}

\begin{definition}[\protect{Nazarov~\cite{Nazarov}};  Ariki, Mathas, Rui ~\cite{ariki-mathas-rui}]
\label{Waff relations}
 The  {\em degenerate affine BMW algebra}
$\Waff_{n, S}=\Waff_{n, S}(\Omega)$ is 
the unital associative $R$--algebra with generators
$\{s_i,e_i, x_j : \ 1\le i<n \text{ and }1\le j\le n \} $
and relations:
\begin{enumerate}
    \item (Involutions)\ \ 
$s_i^2=1$, for $1\le i<n$.
    \item (Affine braid relations)
\begin{enumerate}
\item $s_is_j=s_js_i$ if $|i-j|>1$,
\item $s_is_{i+1}s_i=s_{i+1}s_is_{i+1}$,  for $1\le i<n-1$,
\item $s_ix_j=x_js_i$ if $j\ne i,i+1$.
\end{enumerate}
    \item (Idempotent relations)\ \ 
$e_i^2=\omega_0e_i$, \ \  for $1\le i<n$.
    \item (Commutation relations)
\begin{enumerate}
\item $s_ie_j=e_js_i$, if $|i-j|>1$,
\item $e_ie_j=e_je_i$, if $|i-j|>1$,
\item $e_ix_j=x_je_i$,  if $j\ne i,i+1$,
\item $x_ix_j=x_jx_i$,  for $1\le i,j\le n$.
\end{enumerate}
    \item (Skein relations)\ \ 
        $s_ix_i-x_{i+1}s_i=e_i-1$,  and\ 
        $x_is_i-s_ix_{i+1}=e_i-1$,\  for 
        $1\le i<n$.
    \item (Unwrapping relations)\ \ 
        $e_1x_1^ae_1=\omega_ae_1$, for $a>0$.
    \item (Tangle relations)
\begin{enumerate}
\item $e_is_i=e_i=s_ie_i$,  for $1\le i\le n-1$,
\item $s_ie_{i+1}e_i=s_{i+1}e_i$, and  $e_ie_{i+1}s_i=e_i s_{i+1}$,  for $1\le i\le n-2$,
\item $e_{i+1}e_is_{i+1}=e_{i+1}s_i$,  and $s_{i+1}e_ie_{i+1}=s_i e_{i+1}$,  for $1\le i\le n-2$.
\end{enumerate}
    \item (Untwisting relations)\ \ 
        $e_{i+1}e_ie_{i+1}=e_{i+1}$,  and 
        $e_ie_{i+1}e_i=e_i$, for  
        $1\le i\le n-2$.
    \item (Anti--symmetry relations)\ \ 
        $e_i(x_i+x_{i+1})=0$,  and         $(x_i+x_{i+1})e_i=0$, for 
        $1\le i<n$.
\end{enumerate}
\end{definition}

\begin{definition}[Ariki, Mathas, Rui ~\cite{ariki-mathas-rui}]  \label{cyclo NZ}
Fix an integer $r\ge1$ and elements $u_1,\dots,u_r $ in $S$.  
The {\em degenerate cyclotomic  BMW  algebra} $\W_{n, S, r}=  \W_{n, S, r}(u_1, \dots, u_r)$ is 
the quotient of the degenerate affine BMW algebra $\Waff_{n, S}(\Omega)$ by the relation
$(x_1-u_1)\dots(x_1-u_r) = 0.$
\end{definition}

Due to the symmetry of the relations,  $\Waff_{n, S}$  has a unique $S$--linear  algebra involution $i$ fixing each of the generators.   The involution passes to cyclotomic quotients.  

\subsubsection{Admissibility}  As for the cyclotomic BMW algebras,  to understand the degenerate cyclotomic BMW algebras it is crucial to first understand the ``optimal" case, namely  the case that  $\W_{r, 2} e_1$  is  free  of rank $r$.   We say that the parameters are {\em admissible} if this condition holds. 

It has been shown in ~\cite{goodman-degenerate-admissibility}  and ~\cite{semi-admissible} that
admissibility is equivalent to certain  polynomial conditions on the parameters that were proposed by
Ariki, Mathas and Rui ~\cite{ariki-mathas-rui}, called $u$--admissibility, and also to 
an analogue of the admissibility condition of Wilcox and Yu ~\cite{Wilcox-Yu}
for the cyclotomic BMW algebras.  
It was  shown  in ~\cite{semi-admissible} that  the structure of degenerate cyclotomic BMW algebras with non--admissible parameters can be derived from the admissible case.

In an admissible ground ring, the parameters $\omega_a$  are given by specific polynomial functions of $u_1, \dots, u_r$.  
There is a generic admissible ground ring $R = \Z[\ubold_1, \dots, \ubold_r]$,  where
the $\ubold_j$ are algebraically independent indeterminants.  

\begin{remark}  In previous work on the degenerate cyclotomic BMW algebras, it was always assumed that $2$ is invertible in the ground ring.  However, it was shown in ~\cite{semi-admissible} that this assumption could be eliminated.  In particular, there is no need to adjoin $1/2$ to the generic admissible ground ring $R$.
\end{remark}

\subsubsection{Some basic properties of degenerate cyclotomic BMW algebras}
We 
 establish some elementary properties of degenerate cyclotomic BMW algebras.  Several of the properties can be shown for degenerate affine BMW algebras instead. 
 Let $S$ be any appropriate ground ring for the degenerate affine or cyclotomic BMW algebras, and write $\Waff_n$  for $\Waff_{n, S}$  and $\W_n$  for $\W_{n, S, r}$.  

\begin{lemma}[see ~\cite{ariki-mathas-rui}, Lemma 2.3]                        \label{SX^a}
In the  affine BMW algebra $\Waff_n$,   for  $1\le i<n$ and  $a\ge1$, one has
\begin{equation} \label{equation: S X^a}
s_ix_i^a=x_{i+1}^as_i+\sum_{b=1}^ax_{i+1}^{b-1}(e_i-1)x_i^{a-b}.
\end{equation}
\end{lemma}

\begin{lemma} \label{lemma:  reduced form of words in deg. affine BMW}
 For $n \ge 1$,  $\Waff_{n}$ is contained in the span of $\Waff_{n-1}$ and of elements
of the form  $a \chi_n b$,  where $a, b \in \Waff_{n-1}$  and $\chi_n \in \{e_{n-1}, s_{n-1}, x_n^\alpha :  \alpha \ge 1\}$.  
\end{lemma}

\begin{proof}  We do this by induction on $n$. The base case $n = 1$ is clear since
$\W_{1, S, r}$  is generated by $x_1$.    Suppose now that $n > 1$ and make the appropriate induction hypothesis.    We have to show that a word in the generators having at least
two occurrences of  $e_{n-1}, s_{n-1}$, or a power of $x_n$ can be rewritten as a linear combination of words with fewer occurrences.   

Consider a subword $\chi_n y \chi_n'$  with
$\chi_n, \chi_n'  \in \{e_{n-1}, s_{n-1}, x_n^\alpha :  \alpha \ge 1\}$ and $y \in \Waff_{n-1}$.
If one of $\chi_n, \chi_n'$ is a power of $x_n$, then it commutes with $y$;  say without loss of generality $\chi_n = x_n^\alpha$.  Then  $\chi_n y \chi_n' = y  x_n^\alpha \chi_n'$.
Now consider the cases $\chi_n' = e_{n-1}$,  $\chi_n' = s_{n-1}$,  and $\chi_n' = x_n^\beta$. 
We have $y x_n^\alpha e_{n-1} = y (-1)^\alpha x_{n-1}^\alpha e_{n-1}$  and $y x_n^\alpha x_n^\beta = y x_n^{\alpha + \beta}$.  Finally $y x_n^\alpha s_{n-1}$  can be dealt with using
Lemma \ref{SX^a}.

Suppose both of $\chi_n, \chi_n'$ are in $\{e_{n-1}, s_{n-1}\}$.  If $y \in \Waff_{n-2, S}$, then
 $\chi_n y \chi_n' = y \chi_n  \chi_n'$.  But the product of any two of $e_{n-1}, s_{n-1}$ is either
 $1$ or a multiple of $e_{n-1}$.   If $y \not\in  \Waff_{n-2, S}$, then we can assume, using the induction hypothesis, that $y = y'  \chi y''$,   where $y', y'' \in  \Waff_{n-2, S}$, and
 $\chi$ is one of $e_{n-2}, s_{n-2}$,  or a power of $x_{n-1}$.  Since $\chi_n, \chi_n'$ commute with $y', y''$,  we are reduced to considering $\chi_n \chi \chi_n'$.   Moreover, if $\chi$ is not
 a power of $x_{n-1}$, then essentially we are dealing with a computation in the Brauer algebra, which was done in  ~\cite{Wenzl-Brauer}, Proposition 2.1.  If one of $\chi_n, \chi_n'$ is $s_{n-1}$,  then the computation can be done using Lemma \ref{SX^a}.
 Thus the only interesting
 case is  $e_{n-1} x_{n-1}^\alpha e_{n-1}$.  But by Lemma 4.15 in ~\cite{ariki-mathas-rui},   
 $e_{n-1} x_{n-1}^\alpha e_{n-1} = \omega e_{n-1}$, where $\omega$ is in the center of
 $\Waff_{n-2}$.
 \end{proof}

\begin{lemma} \label{lemma:  deg cyclotomic BMW  axiom 6} 
\mbox{}
\begin{enumerate}
\item 
For $n \ge 3$, $e_{n-1} \Waff_{n-1} e_{n-1} = \Waff_{n-2} e_{n-1}$.
\item  $e_1 \Waff_{1}  e_1 =  \langle \omega_j: j \ge 0 \rangle  \  e_1$, where
$ \langle \omega_j: j \ge 0 \rangle $ denotes the ideal in $S$ generated by all $\omega_j$.  
\item  For $n \ge 2$,  
$e_{n-1}$  commutes with $ \Waff_{n-2} $.   
\end{enumerate}
\end{lemma}

\begin{proof}  First we have to show that if $y \in \Waff_{n-1, S}$, then 
$e_{n-1} y e_{n-1} \in \Waff_{n-2, S} e_{n-1}$.    Using Lemma \ref{lemma:  reduced form of words in deg. affine BMW}, we can suppose that either $y \in \Waff_{n-2, S}$  or 
$y = y' \chi_{n-1} y''$, with $y', y'' \in \Waff_{n-2, S}$,  and $\chi_{n-1} \in \{ e_{n-2}, s_{n-2},  x_{n-1}^\alpha : \alpha \ge 1\}$.   For $\chi_{n-1}$ a power of $x_{n-1}$,   apply Lemma 4.15 from  ~\cite{ariki-mathas-rui}.  In all other cases, the result follows from the defining relations
of $\Waff_n$.   Thus we have $e_{n-1} \Waff_{n-1, S} e_{n-1} \subseteq \Waff_{n-2, S} e_{n-1}$.  For the opposite inclusion, let $x \in \Waff_{n-2, S}$.  Then
$x e_{n-1} =  e_{n-1} x e_{n-2} e_{n-1} \in e_{n-1} \Waff_{n-1, S} e_{n-1}$.  Points (2) and (3) are obvious.
\end{proof}

\begin{lemma}  \label{lemma: Wn e n-1 = W n-1 e n-1}

For $n\ge 2$,   
 $\W_n e_{n-1} = \W_{n-1} e_{n-1}$.  
  \end{lemma}

\begin{proof}  The proof  is similar to the proof of Lemma 5.3 in ~\cite{GG1}.  Using Lemma \ref{lemma:  reduced form of words in deg. affine BMW},  if $x \in \W_n$ and
$x \not\in \W_{n-1}$,  then we can assume that $x = y' \chi_n y''$,  with $y',y'' \in \W_{n-1}$, 
and $\chi_n \in \{e_{n-1}, s_{n-1},  x_n^\alpha: \alpha \ge 1\}$.  Likewise,  we can assume that either $y'' \in \W_{n-2}$ or that $y'' = z' \chi_{n-1} z''$ with $z', z'' \in \W_{n-2}$ and
$\chi_{n-1} \in  \{e_{n-2}, s_{n-2},  x_{n-1}^\beta: \beta \ge 1\}$. The problem reduces to 
showing that $\chi_n e_{n-1}$  and $\chi_{n} \chi_{n-1} e_{n-1}$ lie in $\W_{n-1} e_{n-1}$ for the various choices of $\chi_n, \chi_{n-1}$.   Most of the cases follow directly from the defining relations, while $s_{n-1} x_{n-1}^\beta e_{n-1}$ must be reduced using Lemma \ref{SX^a}, 
and $e_{n-1} x_{n-1}^\beta e_{n-1}$ requires the use of Lemma 4.15 in ~\cite{ariki-mathas-rui}. 
\end{proof} 

\begin{lemma} \label{lemma: x to x e n injective for deg cyclotomic BMW}
 Let $R$ be the universal admissible ring.  For $n \ge 1$,  the  map $x  \mapsto 
x e_{n}$  is injective from $\W_{n, R, r}$  to $\W_{n, R, r} e_{n}$.  
\end{lemma}

\begin{proof} Note that $e_{n+1} (x e_{n}) e_{n+1} = x e_{n+1}$, so it suffices to show that $x \mapsto x e_{n+1}$ is injective.  It follows from Proposition 2.15 and Theorem A in 
~\cite{ariki-mathas-rui} that $\W_{n, R, r}$ has a basis of ``$r$--regular monomials".
The map $x \mapsto x e_{n+1}$  takes the basis elements of $\W_{n, R, r}$ to distinct basis
elements of $\W_{n+2, R, r}$, so is injective. 
\end{proof}

\subsubsection{Degenerate cyclotomic  Hecke algebras}  

\begin{definition}
\label{ahec relations}  Let $S$ be a commutative ring with identity.
 The  {\em degenerate affine  Hecke algebra}
$\ahec {n, S}$ is 
the unital associative $S$--algebra with generators
$$\{s_i,  x_j : \ 1\le i<n \text{ and }1\le j\le n \}, $$
and relations:
\begin{enumerate}
    \item (Involutions)\ \ 
$s_i^2=1$, for $1\le i<n$.
    \item (Affine braid relations)
\begin{enumerate}
\item $s_is_j=s_js_i$ if $|i-j|>1$,
\item $s_is_{i+1}s_i=s_{i+1}s_is_{i+1}$,  for $1\le i<n-1$,
\end{enumerate}
    \item (Commutation relations)
$x_ix_j=x_jx_i$,  for $1\le i,j\le n$ and $s_ix_j=x_js_i$ if $j\ne i,i+1$.
    \item (Skein relations)\ \ 
        $s_ix_i-x_{i+1}s_i=-1$,  and\ 
        $x_is_i-s_ix_{i+1}=-1$,\  for 
        $1\le i<n$.
  \end{enumerate}
  Let $u_1, \dots, u_r$ be elements of $S$.    The {\em degenerate cyclotomic Hecke algebra}
   $\hec {n, S, r}(u_1, \dots, u_r)$ is the quotient of $\ahec n$  by the relation $(x_1 - u_1)(x_2 - u_2)\cdots(x_1 - u_r) = 0.$
   \end{definition}
   
   The degenerate cyclotomic Hecke algebra is a free $S$--module of rank $r^n n!$,  and
   $\hec {n, S, r}(u_1, \dots, u_r)  \break \hookrightarrow \hec {n+1, S, r}(u_1, \dots, u_r)$ for all $n$
   ~\cite{kleshchev-book}.    $\hec {n, S, r}(u_1, \dots, u_r)$ has a unique algebra involution $i$ fixing the generators;  the involutions on the tower of degenerate cyclotomic Hecke algebras are consistent.
   
   It is observed in ~\cite{ariki-mathas-rui}, Section 6,  that the Murphy type cellular basis of the cyclotomic Hecke algebra from ~\cite{dipper-james-mathas} can be easily adapted to the degenerate cyclotomic Hecke algebras.  
Recall that the partially ordered set $\La_n \spp 0$
in the cell datum for 
$\hec{n, S, r} = \hec{n, S, r}(q; u_1, \dots, u_r)$  is the set of $r$--tuples of Young diagrams with total size $n$, ordered by dominance.   For each $\labold \in \La_n \spp 0$,   the
index set $\mathcal T(\labold)$ in the cell datum is the set of standard tableaux of shape
$\labold$.  The proof of strong coherence of the sequence of cyclotomic Hecke algebras
in ~\cite{ariki-mathas},  Proposition 1.9, and ~\cite{mathas-2009}  also applies to the degenerate cyclotomic Hecke algebras.  
   
  Let $J$ be the ideal in the degenerate cyclotomic BMW algebra $\W_{n, S, r}(u_1, \dots, u_r)$ generated by $e_{n-1}$.   It is straightforward to  show that $\W_{n, S, r}(u_1, \dots, u_r)/J \cong
  \hec {n, S, r}(u_1, \dots, u_r)$, as algebras with involution. 
  
  \subsubsection{Verification of the framework axioms for the degenerate cyclotomic BMW algebras}
  
  Let $R$ be the generic admissible integral ground ring, $R =  \Z[ \ubold_1, \dots, \ubold_r]$.    In this section, 
we write  $\W_n$  for $\W_{n, R, r}(\ubold_1, \dots,  \ubold_r)$ and $\hec n$  for
 $\hec{n, R, r}( \ubold_1, \dots, \ubold_r)$.  The field of fractions of $R$ is
$F = \Q( \ubold_1, \dots, \ubold_r)$.    

\begin{proposition}\label{proposition: framework axioms for BMW}
 The two sequences of algebras $(\bmw n)_{n \ge 0}$  and $(H_n)_{n \ge 0}$ satisfy the strong framework axioms of Section \ref{subsection: framework axioms}.
\end{proposition}

\begin{proof}  As observed above,   $(H_n)_{n \ge 0}$ is a strongly coherent tower of cellular algebras, so the strong version of axiom (\ref{axiom Hn coherent}) holds.  Axioms (\ref{axiom: involution on An}) and (\ref{axiom: A0 and A1})  are evident.
$\bmw n^F$ is semisimple by  ~\cite{ariki-mathas-rui}, Theorem 5.3.  Thus axiom (\ref{axiom: semisimplicity})  holds.

We observed above that $\bmw n/\bmw n e_{n-1} \bmw n \cong H_n$, as algebras with involutions.  Thus axiom (\ref{axiom:  idempotent and Hn as quotient of An}) holds.  Axiom 
(\ref{axiom: en An en})  follows from  Lemma \ref{lemma:  deg cyclotomic BMW  axiom 6} and axiom (\ref{axiom:  An en}) from Lemmas \ref{lemma: Wn e n-1 = W n-1 e n-1} and 
\ref {lemma: x to x e n injective for deg cyclotomic BMW}.
Finally, axiom (\ref{axiom: e(n-1) in An en An}) holds  because of the relation $e_{n-1} e_n e_{n-1} = e_{n-1}$.
\end{proof}

\begin{corollary}  Let $S$ be any admissible  ground ring.
The sequence of  degenerate cyclotomic BMW algebras $(\W_{n, S, r})_{n \ge 0}$ is a strongly coherent tower of cellular algebras.
$\W_{n, S, r}$  has cell modules indexed by all  pairs $(\labold,n)$, where $\labold$ is an $r$--tuple of Young diagrams of total size $n$, $n-2$, $n-4, \dots$.   The cell module labeled by
 $(\labold, n)$ has a basis labeled by up--down tableaux of length $n$ and shape $\labold$.
\end{corollary}

Cellularity of degenerate cyclotomic BMW algebras was proved in ~\cite{ariki-mathas-rui}, Section 7.   The cell filtration for restricted modules was proved in ~\cite{rui-si-degnerate}, Theorem 4.15.
The proof of both results here is shorter.  

 \subsubsection{JM elements for degenerate cyclotomic BMW and Hecke algebras}
  The analogue of Jucys--Murphy elements for the degenerate cyclotomic Hecke algebras
  $\hec {n, S, r} =   \hec {n, S, r}(u_1, \dots, u_r)$   are just the generators $x_k$.  In order to eventually distinguish between JM elements in the 
  degenerate cyclotomic Hecke algebras and the degenerate cyclotomic BMW algebras, let us introduce the slightly superfluous notation $L_j\spp 0 = x_j$.   It follows from the defining relations that $L_1 \spp 0 + \dots + L_n\spp 0$ is central in  $\hec {n, S, r}$.  

For an $r$--tuple of Young diagrams
$\labold$ of total size $n$ and a cell $x \in \labold$, the additive  content of the cell is
$$\kappa(x) = u_j  + b - a$$  if $x$ is in row $a$ and column $b$ of the $j$--th component of $\labold$.   For a standard tableau $\mathfrak t$  of shape $\labold$,  and $1 \le j \le n$,  let  $\kappa(j, \mathfrak t) =
\kappa(x)$, where $x$ is the cell occupied by $j$ in $\mathfrak t$.   Let $\{ a_{\mathfrak t}^\labold\}$  be the Murphy type basis of the cell module $\Delta^\labold$  indexed by standard tableaux of shape $\labold$.   Then $L_j\spp 0$ acts by
 \begin{equation} \label{equation: triangular action of Lj cyclotomic Hecke case}
   L_j \spp 0a_\mathfrak t^\labold = \kappa(j, \mathfrak t)\  a_\mathfrak t^\labold + \sum_{\mathfrak s\, \rhd \, \mathfrak t} r_\mathfrak s a_\mathfrak s^\labold,
 \end{equation}
 where the sum is over standard tableaux $\mathfrak s$ greater than $\mathfrak t$ in dominance order  (hence in lexicographic order).   It is noted in ~\cite{ariki-mathas-rui}, Lemma 6.6, that this follows by the argument in
 ~\cite{jantzen-sum-formula}, Section 3.   It follows that the sum
 $L_1\spp 0 + \cdots + L_n\spp 0$  acts as the scalar $\alpha(\labold) = \sum_{x \in \labold}  \kappa(x)$ on
 the cell module $\Delta^\labold$.  Thus  $\{L_n\spp 0 : n \ge 0\}$ is an additive  JM--family in the strongly coherent tower of cellular algebras   $(\hec {n, S, r})_{n \ge 0}$.

In the degenerate cyclotomic BMW algebras 
 $\W_n = \W_{n, R, r}(\ubold_1, \dots, \ubold_r)$ over the generic integral admissible ground ring $R$,  we define $L_j = x_j$  for $1 \le j \le n$.  
    We have $L_n \in \W_n$  and
 $L_n$ commutes with $\W_{n-1}$.  We have $(L_j+ L_{j+1}) e_j = e_j( L_j +L_{j+1}) = 0$
 by the defining relations.   It is clear that $\pi_n(L_j) = L_j\spp 0$,  where 
 $\pi_n : \W_n \to \hec n = \hec {n, R, r}(\ubold_1, \dots, \ubold_r)$ is the quotient map.
 
 It now follows from Theorem \ref{theorem:  JM elements in basic construction algebras 2} that
$\{L_j : j \ge 0\}$  is an additive JM--family in $(\W_n)_{n\ge 0}$, with the sum
$L_1 + \dots + L_n$  acting by $$ \beta((\labold, n)) :=  \alpha(\labold)$$  on the cell module $\Delta^{(\labold, n)}$, if
$\labold$ is an $r$--tuple of  Young diagrams  of total size $k$.  By Proposition \ref{proposition: triangularity property of  additive JM elements},  the action of the elements $L_j$ on the basis of $\Delta^{(\labold, n)}$ labelled by up--down tableaux is triangular:
 \begin{equation} \label{equation: triangular action of Lj bmw case}
   L_j a_\mathfrak t^{(\labold, n)} = \kappa(j, \mathfrak t)\  a_\mathfrak t^{(\labold, n)} + \sum_{\mathfrak s\, \succ \, \mathfrak t} r_\mathfrak s a_\mathfrak s^{(\labold, n)},
 \end{equation}
with $\kappa(j, \mathfrak t) = {\beta(\mathfrak t(j))}-{\beta(\mathfrak t(j-1))}$,
   for some elements $r_\mathfrak s \in R$, depending on $j$ and $\mathfrak t$. 
   Moreover, if $\mathfrak t(j) = (\nubold, j)$  and $\mathfrak t(j-1) = (\mubold, j-1)$, then
   $|\nubold| = |\mubold| \pm 1$.  If   $|\nubold| = |\mubold| + 1$ and $\nubold \setminus \mubold = x$, 
   where $x$ is the cell in row $a$ and column $b$ of the $\ell$--th  component of $\nubold$, 
    then
   $$
   \kappa(j, \mathfrak t)  = 
  {\alpha(\nubold)}- {\alpha(\mubold)} =  \kappa(x) = \ubold_\ell  + (b -a).
     $$
    If   $|\nubold| = |\mubold| - 1$ and $\mubold \setminus \nubold = x$,  then
   $$
   \kappa(j, \mathfrak t) =
{\alpha(\nubold)}- {\alpha(\mubold)} = - \kappa(x)\inv =  -\ubold_\ell  - (b -a).
     $$
This recovers Theorem 5.12 of Rui and Si ~\cite{rui-si-degnerate}.  

\subsection{The Jones--Temperley--Lieb algebras}

Let $S$ be a commutative ring with identity, with distinguished element $\delta$.  The Jones--Temperley--Lieb algebra $\tl_n(S, \delta)$ is the unital $S$--algebra with generators $e_1, \dots, e_{n-1}$ satisfying the relation:
\begin{enumerate}
\item $e_j^2 = \delta e_j$,
\item $e_j e_{j \pm 1} e_j = e_j$,
\item  $e_j e_k =  e_k e_j$, if $|j - k| \ge 2$,
\end{enumerate}
whenever all indices involved are in the range from $1$ to $n-1$.

The Jones--Temperley--Lieb algebra can also be realized as the subalgebra of the Brauer algebra, with
parameter $\delta$,  spanned by  Brauer diagrams {\em without crossings}.   If $J_n$   denotes the ideal in $\tl_n(S, \delta)$ generated by $e_{n-1}$ (or, equivalently, by any $e_j$),  then $\tl_n(S, \delta)/J_n \cong S$.  

The generic ground ring for the Jones--Temperley--Lieb algebras is $R_0 = \Z[\deltabold]$,  where 
$\deltabold$ is an indeterminant over $\Z$.  
 It is shown in ~\cite{GG1}, Section 5.3, that the pair of towers of algebras
 $(\tl_n(R_0, \deltabold))_{n \ge 0}$  and $(R_0)_{n \ge 0}$  satisfies the framework axioms of Section 
\ref{subsection: framework axioms}.  It follows from Theorem \ref{main theorem} that  the sequence of Jones--Temperley--Lieb  algebras is a strongly coherent tower of cellular algebras.   Moreover, the
partially ordered set in the cell datum for $\tl_n$  is naturally realized as 
\begin{equation} \label{equation:  po set for TL cellular structure 1}
\begin{aligned}
&\{ (k, n) : k \le n \text{ and  } 
n - k  \text{ even}\},  \text{ with} \\
 &(k, n) \le  (k', n) \Leftrightarrow k \ge  k'.  
\end{aligned}   
\end{equation}

\begin{proposition} Fix $S$ and $\delta$ and write $\tl_n$  for $\tl_n(S, \delta)$.   For $n \ge 0$  and  $k \le n$,
$\tl_n^{(k, n)}$ is the ideal in $\tl_n$  generated by $e_{k+1} e_{k+3} \cdots e_{n-1}$.  \end{proposition}

\begin{proof}  For $k =n$,  we interpret $e_{k+1} e_{k+3} \cdots e_{n-1}$ as $1$,  so the statement is trivial.   In particular, the statement is true for $n = 0, 1$.   Let $n \ge
 2$ and suppose the statement is
true for $\tl_{n'}$  with $n' < n$.    By the proof of  Theorem 3.2  in ~\cite{GG1},  in particular
Proposition 4.7, for $k < n$  we have  $\tl_n^{(k, n)} = \tl_n e_{n-1}  \tl_{n-2}^{(k, n-2)} \tl_n$.  Applying the
induction hypothesis, 
$$
\begin{aligned}
\tl_n^{(k, n)} &= \tl_n e_{n-1}  \tl_{n-2}^{(k, n-2)} \tl_n \\& = 
\tl_n e_{n-1}  \tl_{n-2} (e_{k+1} e_{k+3} \cdots e_{n-3}) \tl_{n-2} \tl_n  \\ &=
\tl_n  (e_{k+1} e_{k+3} \cdots e_{n-3}e_{n-1})  \tl_n .
\end{aligned}
$$
\end{proof}

\def\powerhalf{^{1/2}}
\def\powerminushalf{^{-1/2}}
\def\pmhalf{^{\pm 1/2}}

Let $R_0$ be as above, and let $\qbold\powerhalf$  be a solution to $\qbold\powerhalf + \qbold\powerminushalf = \deltabold$ in an extension of $R_0$.  Define $R = \Z[\qbold\pmhalf]$ and let $F = \Q(\qbold\pmhalf)$.   Let $H_n$ denote the Hecke algebra $H_{n, R}(\qbold)$.  Then
$\varphi:T_j \mapsto \qbold\powerhalf e_j - 1$  defines a homomorphism from $H_{n, R}(\qbold)$ to
 $\tl_n(R, \deltabold)$, respecting the algebra involutions.     The kernel of $\varphi$ is the ideal in $H_n$  generated by
\begin{equation} \label{equation: ker of map from Hecke to TL}
\xi = T_1 T_2 T_1 + T_1 T_2 + T_2 T_1 + T_1 + T_2 + 1,
\end{equation}
see ~\cite{GHJ}, Corollary 2.11.2. 

Recall from Example \ref{example: hecke algebras 1} that the Hecke algebra $H_n$  has a cell datum  whose  partially ordered set is  the set $Y_n$ of Young diagrams of size  $n$  with
dominance order.  The set $\Gamma_n$  of Young diagrams with at least three columns is an order ideal
in $Y_n$;  let $I_n = H_n(\Gamma_n)$  denote the corresponding $i$--invariant  two sided ideal of $H_n$.

\medskip

The proof of the following lemma is straightforward.

\begin{lemma}  Let  $A$ be a cellular algebra.  Let $\Lambda$  denote the partially ordered set in the cell datum for $A$, let $\Gamma$  be an order ideal in $\Lambda$, and let $A(\Gamma)$  be the corresponding ideal of $A$.   Then $A/A(\Gamma)$ is
a cellular algebra,  with cellular basis $\{c_{\mathfrak s, \mathfrak t}^\la  + A(\Gamma): \la \in  \Lambda \setminus \Gamma; \  \mathfrak s, \mathfrak t \in \mathcal T(\la)\}$.
  \end{lemma}

Applying the lemma to the Hecke algebra, we have that $H_n/I_n$ is a cellular algebra,  with cellular basis 
 $\{m_{\mathfrak s, \mathfrak t}^\la + H_n(\Gamma_n) : \la \in Y_n \setminus \Gamma_n; \    \mathfrak s, \mathfrak t \in \mathcal T(\la)\}$.  The set $Y_n \setminus \Gamma_n$ is the set of Young diagrams of size $n$ with no more than 2 columns. It is totally ordered by dominance.  Write $\la(k, n) = 
 (2^{({n-k})/{2}},  1^{k})$, i.e. the Young diagram with $(n-k)/2$ rows with two boxes and $k$ rows with one box.   Then  
\begin{equation}
\begin{aligned}
&Y_n \setminus \Gamma_n = \{\la(k, n) : k \le n \text{ and } n-k \text{ even}\},  \text{ with} \\
 &\la(k, n) \unlhd \la (k', n) \Leftrightarrow k \ge k';  
\end{aligned}   
\end{equation} 
compare (\ref{equation:  po set for TL cellular structure 1}).

\begin{lemma} $H_n/I_n \cong \tl_n(R, \deltabold)$. 
\end{lemma}

\begin{proof}  For $n = 1, 2$,  $\Gamma_n = \emptyset$ and $I_n  = (0)$.  On the other hand,
$H_n \cong  \tl_n(R, \deltabold) \cong R$.
For $n \ge 3$, let $\mu = (3, 1^{n-3})$. 
 In the notation of ~\cite{Mathas-book}, chapter 3, 
$\xi = m_\mu = m_{\mathfrak t^\mu, \mathfrak t^\mu}^\mu \in  I_n$, where $\xi$ is the element in Equation (\ref {equation: ker of map from Hecke to TL}).     Hence the ideal
$\langle \xi \rangle$  generated  by $\xi$ in $H_n$   is contained in $I_n$.  Therefore, we have 
a surjective homomorphism of involutive algebras  $A_n \cong H_n/\langle \xi \rangle \to H_n/I_n$.   Both algebras are free of rank
$\sum_\la (f_\la)^2 = \frac{1}{n+1} {{2n}\choose{n}}$,  where the sum is over Young diagrams of
size $n$ and no more than two columns, and $f_\la$ is the number of standard Young tableaux of shape $\la$.  Hence, the homomorphism is an isomorphism. 
\end{proof}

We identify $H_n/I_n$ with $A_n$.  By slight abuse of notation, we write $T_j$ for the image of
$T_j$ in $A_n$,  namely $T_j =  \qbold\powerhalf e_j - 1$.  Thus $T_j + 1 = \qbold\powerhalf e_j$. 
We now have potentially two cellular structures on $A_n$,  one inherited from the Hecke algebra and
one obtained by the construction of ~\cite{GG1}, Section 5.3.  

By the description of the cellular structure on the Hecke algebra in ~\cite{Mathas-book}, chapter 3, 
we have that $A_n^{\la(k, n)}$ is the span of $A_n m_{\la(j, n)} A_n$  with
$j \le k$,  where $$m_{\la(j, n)}= (1 + T_1) (1+ T_3) \cdots (1+ T_{n-j-1}) = \qbold^{(n-j)/2} e_1 e_3 \cdots e_{n-j-1}.$$   Thus, in fact,
$$
\begin{aligned}
A_n^{\la(k, n)} &= A_n (e_1 \cdots e_{n-k-1}) A_n \\
&= A_n (e_{k+1} \cdots e_{n-1}) A_n = A_n^{(k, n)}.
\end{aligned}
$$
Moreover, the cell modules from the two cellular structures are explicitly isomorphic:
$$
\begin{aligned}
\Delta^{\lambda(k, n)} &= A_n (e_{1} \cdots e_{n-k-1}) + \breve A_n^{\lambda(k, n)}\\
&\cong  A_n (e_{k+1} \cdots e_{n-1})   + \breve A_n^{(k, n)} = \Delta^{(k, n)}.
\end{aligned}
$$

We can now import the  JM elements from the Hecke algebras (see Example \ref{example: JM elements in Hecke algebra})  to the Jones--Temperley--Lieb algebras. Set 
$L_1 = 1$  and $L_{j+1} =q\inv T_j L_j T_j$   for $j \ge 1$.     Since the cell modules for the  Jones--Temperley--Lieb algebra $A_n$  are in fact cell modules for the Hecke algebra $H_n$, the triangularity property (\ref{equation: triangular action of JM elements in Hecke example})  follows,  and 
 the product $ \prod_{j = 1}^n L_j$  acts as the  scalar $$\alpha(\la(k,n)) = q^{\sum_{x \in \la(k,n)}\kappa(x)}$$  on the cell module $\Delta^{\la(k,n)} = \Delta^{(k,n)}$.
 One can check that $$\frac{\alpha(\la(k,n))}{\alpha(\la(k,n-2))} = q^{-n + 3},$$   independent of $k$, for $n \ge 2$.  It follows from this that $L_n L_{n+1} e_n = e_n L_n L_{n+1} = q^{-n + 2} e_n$  for $n \ge 1$.

\begin{remark}  The same or similar analogues of Jucys--Murphy elements for the Jones--Temperley--Lieb algebras have been considered in ~\cite{ram-halverson-tl} and ~\cite{enyang-tl}.  Those in  ~\cite{enyang-tl} are defined
over the generic ring $R_0 = \Z[\deltabold]$,  but it is not clear that they have, or can be modified to have,  the multiplicative property
(resp. additive property) of Definition \ref{definition: multiplicative JM family}  or \ref{definition: additive JM family}.  
\end{remark}  

\bibliographystyle{amsplain}
\bibliography{JM} 
\end{document}